\documentclass[11pt]{amsart}
\usepackage{amssymb}
\usepackage{graphicx}                            
\usepackage{hyperref}
\usepackage{setspace}                            

\newtheorem{lemma}{Lemma}[section]
\newtheorem{theorem}{Theorem}[section]
\newtheorem{proposition}{Proposition}[section]

\newcommand{\adfQED}{\hfill $\square$}           
\newcommand{\adfmod}[1]{~(\mathrm{mod}~#1)}
\newcommand{\adfhat}[1]{\hat{#1}}
\newcommand{\adfhide}[1]{}
\newcommand{\adfunhide}[1]{#1}
\newcommand{\ADFvfyParStart}[1]{}
\newcommand{\adfDgap}{\vskip 2mm}              
\newcommand{\adfLgap}{\vskip 0mm}              
\newcommand{\adfsplit}{\par}                   
\newcommand{\adfBfont}{\normalsize}            

\begin{document}
\title{Design spectra for 6-regular graphs with 12 vertices}
\author{A. D. Forbes}
\address{LSBU Business School,
London South Bank University,
103 Borough Road,
London SE1 0AA, UK.}
\email{anthony.d.forbes@gmail.com}
\author{C. G. Rutherford}
\address{LSBU Business School,
London South Bank University,
103 Borough Road,
London SE1 0AA, UK.}
\email{c.g.rutherford@lsbu.ac.uk}
\date{\today} 

\begin{abstract}
The design spectrum of a simple graph $G$ is the set of positive integers $n$ such that there exists an edgewise
decomposition of the complete graph $K_n$ into $n(n - 1)/(2 |E(G)|)$ copies of $G$.
We compute the design spectra for 7788 6-regular graphs with 12 vertices.
\end{abstract}

\maketitle


\section{Introduction}\label{sec:Introduction}
There are 7849 6-regular graphs with 12 vertices.
Of these, 7848 are available online as complements of the connected 12-vertex, 5-regular graphs constructed by Meringer, \cite{Meringer1999}.
Throughout this paper we refer to them by their positions in Meringer's list.
Thus graph $n$, $1 \le n \le 7848$, denotes the complement of the $n$-th graph in
the list of the edge sets of all connected 12-vertex, 5-regular graphs at \cite{MeringerReggraphs1999}.
To the list we append the complete bipartite graph $K_{6,6}$ as number 7849.

If $F$ and $G$ are simple graphs,
an {\em edgewise decomposition} of $F$ into $G$, which we also refer to as
a {\em $G$-decomposition of $F$}, is a partition $\mathcal{E}$ of the edges of $F$
such that each $E \in \mathcal{E}$ is the edge set of a graph isomorphic to $G$.
If $F$ is the complete graph $K_n$, we usually refer to the decomposition as a {\em $G$-design} of order $n$.
The {\em design spectrum} of $G$ is the set of positive integers $n$ for which a $G$-design of order $n$ exists.
If $G$ is $d$-regular, the necessary conditions for the existence of a $G$-design are
\begin{align}
n &\ge |V(G)| \text{~or~} n = 1, \nonumber \\
n(n - 1) &\equiv 0 \adfmod{2|E(G)|}, \label{eqn:necessary-conditions} \\
   n - 1 &\equiv 0 \adfmod{d}. \nonumber
\end{align}

Given a $d$-regular graph $G$, by a theorem of Wilson, \cite{Wilson1976}, the conditions (\ref{eqn:necessary-conditions})
are sufficient for all sufficiently large $n$, and hence
the determination of $G$'s design spectrum is actually a finite problem.
However, it is usually impossible to resolve all of the cases not covered by `sufficiently large'
whenever $d$ or the chromatic number is large.
Nevertheless, design spectra have been computed for many graphs, including some infinite classes.
For example, from the early history of design theory
we know the spectrum for the complete graph $K_k$ when $k = 2$ (trivial),
$k = 3$ (Kirkman, 1847, \cite{Kirkman1847}) and $k = 4, 5$ (Hanani, 1961, \cite{Hanani1961}),
but there are unresolved cases when $k \ge 6$.
%
The Platonic graphs have also received attention.
Apart from the icosahedron they are either 3- or 4-regular and their design spectra have been completely resolved,
\cite{AdamsBillingtonRodger1994, AdamsBryantForbesGriggs2012, BryantZanatiGardner1994, GriggsResminiRosa1992,
Hanani1961, Kotzig1981, Maheo1980}.
On the other hand, the icosahedron is 5-regular and the partial solution of its design spectrum
leaves 6 unresolved cases, \cite{AdamsBryant1996, ForbesGriggs2012}.
A major barrier seems to be the graph's chromatic number.
In the successful examples mentioned above the chromatic number is at most 5, and
we are not aware of any (6 or more)-chromatic graph for which the design spectrum has been determined.
For a survey of the subject to the year 2008, the reader is referred to \cite{AdamsBryantBuchanan2008}.

We summarize our results.
When $G$ is 6-regular and has 12 vertices
the necessary conditions (\ref{eqn:necessary-conditions}) simplify to
\begin{equation}
\label{eqn:v12-6reg-necessary-conditions}
n \equiv 1 \adfmod{72}.
\end{equation}
In Section~\ref{sec:Graph-designs-general-constructions} we prove that
the condition (\ref{eqn:v12-6reg-necessary-conditions}) is sufficient
(and therefore the design spectrum is determined) for 7788 graphs:
\begin{enumerate}
\item[(i)]the 2-chromatic graph, 7849;
\item[(ii)]all of the forty-nine 3-chromatic graphs,
       1,    2,    3,    4,   17,   18,   20,   22,   23,   24,
     201,  203,  206,  207,  228,  312,  527,  529,  590,  599,
     601,  850, 1106, 1233, 1261, 1698, 1702, 1825, 1835, 1839,
    2040, 2045, 2051, 2053, 2471, 2562, 2563, 2574, 2581, 3179,
    3191, 3193, 3241, 3243, 6383, 6385, 6390, 6397, 6401;
\item[(iii)]6487 of the 6498 4-chromatic graphs;
\item[(iv)]1251 of the 1299 5-chromatic graphs.
\end{enumerate}
Also we note that graph 7849 has already been solved by Rosa, \cite{Rosa1967}, \cite[Theorem 5.3]{AdamsBryantBuchanan2008}.
The 61 graphs where we are not entirely successful are as follows.
\begin{enumerate}
\item[(i)]For eleven 4-chromatic graphs,
      10,   13,   59,  130,  211,  432,  551, 3281, 6729, 7679,
    7743,
and forty-three 5-chromatic graphs,
      16,  163,  424,  635,  659,  670,  671,  687,  692,  701,
     702,  707,  722,  733, 1063, 1438, 3101, 3443, 3447, 4001,
    4069, 4070, 4074, 4096, 4108, 4317, 4764, 4778, 5701, 5859,
    5913, 6339, 6391, 6657, 6751, 7353, 7421, 7531, 7603, 7667,
    7752, 7761, 7803,
the necessary condition (\ref{eqn:v12-6reg-necessary-conditions}) is sufficient with the possible exception of order 505.
\item[(ii)]There are five 5-chromatic graphs where we were unable to obtain a key decomposition, namely that of $K_{18^5}$,
and so we fall somewhat short of obtaining their design spectra:
graphs 672, 716, 6187, 6196, 7824.
\item[(iii)]We made no attempt to address the two 6-chromatic graphs: 703 and 7848.
\end{enumerate}

It is obvious to us that providing detailed proofs for all of the 7788 successful cases
would overload the main part of this paper with an enormous amount of data.
Instead we focus our attention on the eleven graphs that are 3- or 4-chromatic and vertex-transitive,
$$201, 6383, 6397, 6401, 6406, 6408, 6753, 7677, 7754, 7845, 7847,$$
as well as a few graphs that require slightly special treatment,
$$1513, 3470, 6713, 7700, 7840.$$
These sixteen graphs are the subject of Theorem~\ref{thm:spectra-complete-main}.
In Section~\ref{sec:Graphs} we specify the graphs as used for our analysis.
They are illustrated in Figures~\ref{fig:graph-pictures-first}--\ref{fig:graph-pictures-last},
where the positioning of the vertices around the circles has been adjusted
to make the pictures of the vertex-transitive ones look pretty.
The only other vertex-transitive graphs are 7848 and 7849, \cite{McKay1979}.

For all other graphs where we have been successful, the construction details have been placed in the appendix,
which is present only in this paper's preprint at
\begin{center}
\url{https://arxiv.org/abs/blah.blah}.
\end{center}
The purpose of the appendix is to provide material that the interested reader can use to verify our claims with the aid of a computer.
The appendix contains the details for each graph $G$ that satisfies one of these conditions:
\begin{enumerate}
\item[(A)]a $G$-decomposition of $K_{24^4}$ is available, 6311 graphs;
\item[(B)]$G$-decompositions of $K_{18^5}$, $K_{6^7}$ and $K_{9^9}$ are available but not $K_{24^4}$, 1471 graphs;
\item[(C)]$G$-decompositions of $K_{18^5}$ and $K_{9^9}$ are available but not
$K_{24^4}$, $K_{18^4}$, $K_{6^7}$ or $K_{72^7}$, $54$ graphs.
\end{enumerate}
They are the subject of Theorems~\ref{thm:spectra-complete-appendix} and \ref{thm:spectra-missing-505-appendix}.

The proofs in Section~\ref{sec:Graph-designs-general-constructions}
employ a technique of design theory known as Wilson's fundamental construction, \cite{Wilson1972}.
The method uses group divisible designs to build large graph decompositions from small ones.
In Section~\ref{sec:Group-divisible-designs} we give the definition of a group divisible design
that is relevant to our paper, and
in Lemma~\ref{lem-GDD-existence} we collect together known existence results for the types that we require.
The sequence of lemmas in Section~\ref{sec:Graph-designs-direct-constructions} provides
the details of direct constructions for decompositions of certain small complete and complete multipartite graphs
into the 12-vertex 6-regular graphs of Section~\ref{sec:Graphs}.
Our main theorems are in Section~\ref{sec:Graph-designs-general-constructions},
and we finish with some informal remarks in Section~\ref{sec:Concluding-remarks}.


\section{Graphs}\label{sec:Graphs}
Here we specify the graphs that we deal with in Sections~\ref{sec:Graph-designs-direct-constructions} and \ref{sec:Graph-designs-general-constructions}.
A graph is coded as an ordered 11-tuple $(a_1, a_2, \dots, a_{11})$, where
the binary digits of $a_i$ constitute row $i$ of the above-diagonal part of the adjacency matrix, $i = 1, 2, \dots, 11$.
The chromatic number is indicated by $\chi$.
The `L' number refers to the corresponding vertex-transitive graph in McKay's list, \cite{McKay1979}.

\begin{enumerate}
\setlength{\itemsep}{1mm}
\item[{\bf 201}] (63,63,207,207,51,51,12,12,0,0,0),
$\chi=3$, L26, circulant-12-2-3-4.
\item[{\bf 1513}] (63,95,175,243,45,30,17,2,4,0,0),
$\chi=4$.
\item[{\bf 3470}] (63,95,399,179,92,51,28,0,4,0,1),
$\chi=4$.
\item[{\bf 6383}] (63,207,243,252,21,42,5,10,1,2,0),
$\chi=3$, L30, circulant-12-1-2-5.
\item[{\bf 6397}] (63,207,243,252,69,10,21,10,4,0,1),
$\chi=3$, L32, circulant-12-1-4-5.
\item[{\bf 6401}] (63,207,243,252,69,18,26,4,5,0,1),
$\chi=3$, L29, circulant-12-2-4-5.
\item[{\bf 6406}] (63,207,343,91,93,33,30,2,4,0,0),
$\chi=4$, L35, complement of (octahedron $\times$ $K_2$).
\item[{\bf 6408}] (63,207,343,91,93,34,30,4,1,0,0),
$\chi=4$, L31.
\item[{\bf 6713}] (63,207,343,121,122,36,3,6,5,0,0),
$\chi=4$.
\item[{\bf 6753}] (63,207,343,171,53,57,6,10,4,0,0),
$\chi=4$, L37, complement of the icosahedron.
\item[{\bf 7677}] (63,207,371,211,116,36,24,8,5,2,1),
$\chi=4$, L27.
\item[{\bf 7700}] (63,207,371,213,92,48,26,0,6,1,1),
$\chi=4$.
\item[{\bf 7754}] (63,207,371,220,116,40,17,2,6,1,1),
$\chi=4$, L36, circulant-12-3-4-5.
\item[{\bf 7840}] (63,207,497,242,84,40,24,4,3,3,1),
$\chi=4$.
\item[{\bf 7845}] (63,455,504,75,116,12,21,10,2,1,1),
$\chi=4$, L33, line graph of the octahedron.
\item[{\bf 7847}] (63,455,504,195,73,24,30,12,0,3,1),
$\chi=4$, L28, circulant-12-2-3-5.
\end{enumerate}

\section{Group divisible designs}\label{sec:Group-divisible-designs}
For the purpose of this paper, a {\em group divisible design}, $K$-GDD, of type $g_1^{u_1} g_2^{u_2} \dots g_r^{u_r}$
is an ordered triple ($V, \mathcal{G}, \mathcal{B}$)
where
\begin{enumerate}
\item[(i)]$V$ is a set of $u_1 g_1 + u_2 g_2 + \dots + u_r g_r$ {\em points},
\item[(ii)]$\mathcal{G}$ is a partition of $V$ into $u_i$ subsets of size $g_i$, $i = 1, 2, \dots, r$, called \textit{groups}, and
\item[(iii)]$\mathcal{B}$ is a collection of subsets of cardinalities $k \in K$, called \textit{blocks}, which has property that each pair of points from distinct groups occurs in precisely one block but
    a pair of distinct points from the same group does not occur in any block.
\end{enumerate}
We usually refer to a $\{k\}$-GDD as $k$-GDD.
A {\em parallel class} in a group divisible design is a subset of blocks that precisely covers the point set.
A $k$-GDD is called {\em resolvable}, and denoted by $k$-RGDD, if the entire set of blocks
can be partitioned into parallel classes.

Our first lemma asserts the existence of the group divisible designs that we require in Section~\ref{sec:Graph-designs-general-constructions}.
%
\begin{lemma}[\cite{
AbelColbournDinitz2007,
AbelGeGreigLing2009,
BrouwerSchrijverHanani1977,
Ge2007,
GeLing2004c,
GeLing2005,
GeMiao2007,
HananiRayChaudhuriWilson1972,
WeiGe2014}] \label{lem-GDD-existence}~
\begin{enumerate}
\item[(i)]There exists a $4$-$\mathrm{GDD}$ of type $g^u$ if
$u \ge 4$, $g(u - 1) \equiv 0 \adfmod{3}$ and $g^2u(u - 1) \equiv 0 \adfmod{12}$,
except for $(g,u) \in \{(2,4), (6,4)\}$.
\item[(ii)]There exists a $4$-$\mathrm{GDD}$ of type $6^u 3^1$ if $u \ge 4$.
\item[(iii)]There exists a $4$-$\mathrm{GDD}$ of type $3^5 6^1$.
\item[(iv)]There exists a $\{4,5\}$-$\mathrm{GDD}$ of type $4^{3t + 1} m^1$ for $m \ge 0$ and $t \ge m/4$.
\item[(v)]If $g \in \{4, 8\}$, there exist $5$-$\mathrm{GDD}$s of types $g^{5t}$ and $g^{5t + 1}$ for $t \ge 1$.
\item[(vi)]There exists
\begin{enumerate}
\item[]a $5$-$\mathrm{GDD}$ of type $4^{5t} 8^1$ for $t \ge 3$,  
\item[]a $5$-$\mathrm{GDD}$ of type $4^{5t} 12^1$ for $t \ge 2$,
and
\item[]a $5$-$\mathrm{GDD}$ of type $4^{5t} 16^1$ for $t \ge 4$. 
\end{enumerate}
\item[(vii)]There exists a $5$-$\mathrm{GDD}$ of type $12^5 16^1$.
\item[(viii)]There exist $7$-$\mathrm{GDD}$s of types $1^7$ and $12^7$.
\item[(ix)]There exists a $9$-$\mathrm{GDD}$ of type $8^9$.
\end{enumerate}
\end{lemma}
%
\begin{proof}
(i)--(iii)~ See \cite{BrouwerSchrijverHanani1977} or \cite{Ge2007}.

(iv)~ There exists a 4-RGDD of type $4^{3t + 1}$ whenever $t \ge 1$, \cite{HananiRayChaudhuriWilson1972}, see also \cite[Theorem IV.5.44]{GeMiao2007}.
There are $4t$ parallel classes, $P_1$, $P_2$, \dots, $P_{4t}$, say.
If $m = 0$, we are done.
If $1 \le m \le 4t$, we add an extra group $\{x_1, x_2, \dots, x_m\}$ of size $m$ and
we augment each block of $P_i$ with $x_i$, $i = 1$, 2, \dots, $m$.
The result is a $\{4,5\}$-GDD of type $4^{3t + 1} m^1$.

(v)~ See \cite{GeLing2005} or \cite{WeiGe2014} or \cite[Theorem IV.4.16]{Ge2007}.

(vi)~ See \cite{AbelGeGreigLing2009}.
A $(v, \{5, w^*\}, 1)$-PBD is a pairwise balanced design on $v$ points where
one block has size $w$ and all others have size 5.
The blocks have the `balanced' property---each pair of points occurs in precisely one block.
Theorems~1 and 30 of \cite{AbelGeGreigLing2009} together assert the existence of
\begin{enumerate}
\item[]a $(20t + 9, \{5, 9^*\}, 1)$-PBD for $t \ge 3$,
\item[]a $(20t + 13, \{5, 13^*\}, 1)$-PBD for $t \ge 2$ and
\item[]a $(20t + 17, \{5, 17^*\}, 1)$-PBD for $t \ge 4$.
\end{enumerate}
Remove a point from the block of size $w \in \{9, 13, 17\}$.
The resulting blocks of sizes $4$ and $w - 1$ form the groups of a 5-GDD of type $4^{5t} (w - 1)^1$.

(vii)~ See \cite{GeLing2004c} 
or \cite[Theorem IV.4.17]{Ge2007}.

(viii)~ The 7-GDD of type $12^7$ is constructed from 5 mutually orthogonal Latin squares of side 12;
see \cite[Table III.3.87]{AbelColbournDinitz2007}.
The other one is trivial.

(ix)~ This follows from the existence of a projective plane of order 8.
%
\end{proof}


\section{Graph designs: direct constructions}\label{sec:Graph-designs-direct-constructions}
In a sequence of lemmas we give the direct constructions of $G$-decompositions
that we require for our proofs in Section~\ref{sec:Graph-designs-general-constructions}.

For graph $n$, the set of labelled graphs that form the decomposition is generated from one or two base blocks by a specified mapping.
A base block is a subscripted ordered 12-tuple
$(\ell_1, \ell_2, \dots, \ell_{12})_n$
where, for $i \in \{1, 2, \dots, 12\}$,
label $\ell_i$ is attached to vertex $i$ of
graph $n$ as defined in Section~\ref{sec:Graphs} or in the appendix.

The expression $x$ mod $y$ is the integer in $\{0, 1, \dots, y-1\}$ that is congruent to $x$ modulo $y$.
%
%
\begin{lemma} \label{lem-Design-24-4}
For each of the 12-vertex, 6-regular graphs
$$201, 6383, 6397, 6401, 6753, 7677, 7754,$$
there exists an edgewise decomposition of the complete multipartite graph $K_{24^4}$ into $96$ copies of the graph.
\end{lemma}
%
\begin{proof}
The point set is $\mathbb{Z}_{96}$ partitioned by residue class modulo 4.
The base blocks are developed by
$x \mapsto x + d \mathrm{~mod~} 96$,
$0 \le d < 96$.


$(17, 33, 57, 41, 63, 39, 6, 54, 3, 67, 10, 58)_{201}$

$(35, 24, 95, 39, 9, 10, 4, 82, 81, 42, 61, 26)_{6383}$

$(53, 69, 61, 76, 59, 10, 67, 54, 15, 26, 35, 38)_{6397}$

$(31, 8, 80, 4, 85, 71, 10, 69, 65, 82, 81, 38)_{6401}$

$(91, 44, 29, 59, 57, 43, 68, 0, 17, 50, 10, 94)_{6753}$

$(3, 91, 20, 27, 40, 10, 1, 74, 36, 9, 94, 73)_{7677}$

$(38, 86, 78, 19, 57, 0, 27, 92, 93, 40, 87, 53)_{7754}$
\ADFvfyParStart{\{\{24, 4\}, \{96, \{\{1, 96, 1, 1, \{\{96, 1\}\}\}\}, \{\{24, 4\}\}\}\}}
\end{proof}
%
\begin{lemma} \label{lem-Design-18-5}
For each of the 12-vertex, 6-regular graphs
$$1513, 3470, 6406, 6408, 6713, 7700, 7840, 7845, 7847,$$
there exists an edgewise decomposition of $K_{18^5}$ into $90$ copies of the graph.
\end{lemma}
%
\begin{proof}
The point set is $\mathbb{Z}_{90}$ partitioned by residue class modulo 5.
The base blocks are developed by
$x \mapsto x + d \mathrm{~mod~} 90$,
$0 \le d < 90$.


$(22, 82, 27, 57, 86, 39, 50, 38, 0, 48, 25, 33)_{1513}$

$(27, 41, 81, 74, 57, 84, 71, 63, 23, 80, 43, 55)_{3470}$

$(19, 45, 9, 35, 22, 6, 25, 26, 88, 11, 57, 63)_{6406}$

$(14, 9, 78, 29, 47, 21, 7, 5, 76, 41, 10, 60)_{6408}$

$(72, 56, 54, 52, 47, 48, 35, 36, 14, 75, 13, 28)_{6713}$

$(2, 10, 81, 25, 79, 72, 15, 73, 61, 63, 78, 14)_{7700}$

$(60, 0, 69, 46, 3, 17, 33, 67, 79, 9, 42, 1)_{7840}$

$(22, 27, 7, 20, 61, 78, 4, 30, 31, 0, 68, 6)_{7845}$

$(3, 53, 58, 72, 25, 74, 34, 40, 47, 2, 11, 89)_{7847}$
\ADFvfyParStart{\{\{18, 5\}, \{90, \{\{1, 90, 1, 1, \{\{90, 1\}\}\}\}, \{\{18, 5\}\}\}\}}
\end{proof}
%
\begin{lemma} \label{lem-Design-6-7}
For each of the 12-vertex, 6-regular graphs
$$6406, 6408, 7845,$$
there exists an edgewise decomposition of $K_{6^7}$ into $21$ copies of the graph.
\end{lemma}
%
\begin{proof}
The point set is $\mathbb{Z}_{42}$ partitioned by residue class modulo 7.
The base blocks are developed by
$x \mapsto x + 2d \mathrm{~mod~} 42$,
$0 \le d < 21$.


$(0, 2, 34, 39, 1, 35, 17, 37, 27, 4, 8, 24)_{6406}$

$(41, 40, 20, 17, 27, 14, 26, 21, 8, 9, 18, 25)_{6408}$

$(26, 19, 27, 7, 24, 25, 18, 28, 10, 9, 15, 38)_{7845}$
\ADFvfyParStart{\{\{6, 7\}, \{42, \{\{1, 21, 1, 1, \{\{42, 2\}\}\}\}, \{\{6, 7\}\}\}\}}
\end{proof}
%
\begin{lemma} \label{lem-Design-9-9}
For each of the 12-vertex, 6-regular graphs
$$1513, 3470, 6406, 6408, 6713, 7700, 7840, 7845, 7847,$$
there exists an edgewise decomposition of $K_{9^9}$ into $81$ copies of the graph.
\end{lemma}
%
\begin{proof}
The point set is $\mathbb{Z}_{81}$ partitioned by residue class modulo 9.
The base blocks are developed by
$x \mapsto x + d \mathrm{~mod~} 81$,
$0 \le d < 81$.


$(36, 37, 76, 45, 77, 20, 64, 16, 6, 79, 60, 71)_{1513}$

$(72, 49, 40, 21, 17, 47, 5, 55, 12, 37, 62, 73)_{3470}$

$(58, 20, 27, 59, 23, 75, 11, 51, 79, 21, 16, 8)_{6406}$

$(60, 53, 8, 1, 28, 58, 38, 72, 12, 57, 0, 14)_{6408}$

$(76, 2, 20, 30, 60, 13, 0, 32, 44, 75, 35, 24)_{6713}$

$(21, 6, 12, 76, 32, 65, 61, 9, 22, 52, 27, 14)_{7700}$

$(0, 62, 31, 79, 12, 23, 35, 15, 22, 5, 28, 60)_{7840}$

$(80, 10, 70, 33, 29, 23, 68, 37, 48, 3, 34, 63)_{7845}$

$(32, 23, 78, 71, 58, 1, 26, 79, 63, 40, 0, 2)_{7847}$
\ADFvfyParStart{\{\{9, 9\}, \{81, \{\{1, 81, 1, 1, \{\{81, 1\}\}\}\}, \{\{9, 9\}\}\}\}}
\end{proof}
%
\begin{lemma} \label{lem-Design-18-4}
For each of the 12-vertex, 6-regular graphs
$$1513, 3470, 6713,$$
there exists an edgewise decomposition of $K_{18^4}$ into $54$ copies of the graph.
\end{lemma}
%
\begin{proof}
The point set is $\mathbb{Z}_{72}$ partitioned by residue class modulo 3 for points $\{0, 1, \dots, 53\}$, and $\{54, 55, \dots, 71\}$.
The base blocks are developed by
$x \mapsto x + d \mathrm{~mod~} 54$ for $0 \le x < 54$,
$x \mapsto (x + d \mathrm{~mod~} 18) + 54$ for $54 \le x < 72$,
$0 \le d < 54$.


$(6, 54, 58, 60, 48, 33, 16, 20, 5, 25, 4, 53)_{1513}$

$(29, 26, 50, 6, 64, 9, 68, 4, 3, 37, 10, 62)_{3470}$

$(68, 56, 54, 49, 22, 8, 6, 3, 0, 37, 29, 50)_{6713}$
\ADFvfyParStart{\{\{18, 4\}, \{72, \{\{1, 54, 1, 1, \{\{54, 1\}, \{18, 1\}\}\}\}, \{\{18, 3\}, \{18, 1\}\}\}\}}
\end{proof}
%
\begin{lemma} \label{lem-Design-72-7}
For each of the 12-vertex, 6-regular graphs
$$7700, 7840, 7847,$$
there exists an edgewise decomposition of $K_{72^7}$ into $3024$ copies of the graph.
\end{lemma}
%
\begin{proof}
The point set is $\mathbb{Z}_{504}$ partitioned by residue class modulo 7.
There are two base blocks for each graph.
They are developed by
$x \mapsto 25^e x + d \mathrm{~mod~} 504$, $0 \le e < 3$, $0 \le d < 504$.


$(341, 413, 36, 142, 235, 339, 156, 111, 99, 335, 270, 178)_{7700}$

$(358, 386, 268, 470, 154, 263, 231, 166, 73, 333, 137, 412)_{7700}$

$(292, 246, 10, 196, 453, 261, 167, 130, 343, 364, 160, 452)_{7840}$

$(131, 29, 120, 150, 226, 385, 130, 363, 318, 122, 172, 128)_{7840}$

$(368, 0, 249, 149, 281, 306, 335, 225, 262, 338, 320, 52)_{7847}$

$(346, 206, 111, 393, 285, 478, 151, 362, 0, 218, 268, 83)_{7847}$
\ADFvfyParStart{\{\{72, 7\}, \{504, \{\{2, 504, 3, 25, \{\{504, 1\}\}\}\}, \{\{72, 7\}\}\}\}}
\end{proof}
%
\begin{lemma} \label{lem-Design-73-et-al}
For each of the 12-vertex, 6-regular graphs
\begin{align*}
201, 1513, 3470, 6383, 6397, 6401, 6406, 6408, \\ 6713, 6753, 7677, 7700, 7754, 7840, 7845, 7847,
\end{align*}
there exist designs of orders
$73$, $145$, $217$, $289$, $433$, $577$ and $1009$.
\end{lemma}
%
\begin{proof}
For a design of order $n$, the point set is $\mathbb{Z}_{n}$.
The blocks are developed from a single base block by
$x \mapsto \omega^e x + d  \mathrm{~mod~} n$, $0 \le  e < (n - 1)/72$, $0 \le d < n$,
where $\omega$ is a specified parameter.
%
%


{Order 73}, $\omega = 1$:


$(0, 1, 2, 3, 4, 23, 32, 67, 40, 62, 19, 26)_{201}$

$(0, 1, 2, 3, 4, 50, 25, 31, 41, 70, 58, 63)_{1513}$

$(0, 1, 2, 3, 12, 4, 47, 49, 62, 69, 21, 35)_{3470}$

$(0, 1, 2, 3, 4, 24, 18, 68, 20, 63, 36, 44)_{6383}$

$(0, 1, 2, 3, 4, 70, 28, 41, 55, 65, 24, 44)_{6397}$

$(0, 1, 2, 3, 4, 15, 7, 51, 41, 20, 65, 47)_{6401}$

$(0, 1, 2, 3, 4, 8, 10, 33, 28, 36, 22, 17)_{6406}$

$(0, 1, 2, 3, 4, 19, 10, 41, 54, 8, 62, 49)_{6408}$

$(0, 1, 2, 3, 4, 18, 69, 52, 46, 20, 64, 40)_{6713}$

$(0, 1, 2, 3, 5, 9, 44, 45, 51, 54, 18, 14)_{6753}$

$(0, 1, 2, 3, 6, 58, 13, 67, 23, 43, 41, 27)_{7677}$

$(0, 1, 2, 3, 5, 10, 53, 26, 11, 40, 56, 15)_{7700}$

$(0, 1, 2, 3, 5, 12, 27, 70, 50, 34, 14, 55)_{7754}$

$(0, 1, 2, 3, 5, 10, 37, 63, 53, 19, 25, 42)_{7840}$

$(0, 1, 2, 3, 5, 54, 32, 63, 19, 60, 26, 65)_{7845}$

$(0, 1, 2, 3, 7, 30, 26, 12, 37, 9, 42, 22)_{7847}$
\ADFvfyParStart{\{73, \{73, \{\{1,73,1,1, \{\{73,1\}\}\}\}\}\}}

{Order 145}, $\omega = 12$:


$(0, 1, 2, 3, 4, 6, 19, 60, 125, 139, 104, 117)_{201}$

$(0, 1, 2, 3, 4, 12, 8, 31, 22, 88, 112, 55)_{1513}$

$(0, 1, 2, 3, 5, 12, 97, 28, 31, 40, 20, 93)_{3470}$

$(0, 1, 2, 3, 4, 8, 106, 48, 23, 68, 19, 115)_{6383}$

$(0, 1, 2, 3, 4, 8, 132, 114, 65, 95, 124, 46)_{6397}$

$(0, 1, 2, 3, 4, 9, 41, 117, 18, 98, 54, 120)_{6401}$

$(0, 1, 2, 3, 4, 9, 110, 34, 81, 129, 123, 108)_{6406}$

$(0, 1, 2, 3, 4, 8, 76, 121, 53, 108, 18, 83)_{6408}$

$(0, 1, 2, 3, 4, 8, 120, 132, 24, 42, 83, 94)_{6713}$

$(0, 1, 2, 3, 5, 7, 87, 137, 37, 104, 54, 23)_{6753}$

$(0, 1, 2, 3, 5, 8, 67, 125, 33, 137, 102, 71)_{7677}$

$(0, 1, 2, 3, 5, 8, 47, 130, 137, 94, 116, 28)_{7700}$

$(0, 1, 2, 3, 5, 10, 21, 31, 105, 70, 34, 93)_{7754}$

$(0, 1, 2, 3, 5, 10, 72, 130, 26, 107, 66, 53)_{7840}$

$(0, 1, 2, 3, 5, 18, 10, 82, 40, 89, 102, 116)_{7845}$

$(0, 1, 2, 3, 6, 12, 131, 64, 86, 96, 42, 117)_{7847}$
\ADFvfyParStart{\{145, \{145, \{\{1, 145, 2, 12, \{\{145, 1\}\}\}\}\}\}}

{Order 217}, $\omega = 25$:


$(0, 1, 2, 3, 4, 6, 115, 206, 157, 196, 40, 90)_{201}$

$(0, 1, 2, 3, 4, 6, 49, 147, 157, 208, 118, 183)_{1513}$

$(0, 1, 2, 3, 5, 9, 40, 24, 12, 120, 147, 152)_{3470}$

$(0, 1, 2, 3, 4, 7, 20, 58, 30, 127, 43, 136)_{6383}$

$(0, 1, 2, 3, 4, 8, 112, 14, 172, 201, 38, 81)_{6397}$

$(0, 1, 2, 3, 4, 8, 40, 124, 154, 47, 19, 55)_{6401}$

$(0, 1, 2, 3, 4, 8, 40, 123, 199, 48, 73, 31)_{6406}$

$(0, 1, 2, 3, 4, 8, 21, 48, 201, 188, 149, 138)_{6408}$

$(0, 1, 2, 3, 4, 8, 101, 15, 147, 168, 47, 200)_{6713}$

$(0, 1, 2, 3, 5, 7, 14, 123, 186, 198, 172, 146)_{6753}$

$(0, 1, 2, 3, 5, 8, 38, 166, 107, 85, 21, 124)_{7677}$

$(0, 1, 2, 3, 5, 8, 47, 65, 136, 33, 205, 196)_{7700}$

$(0, 1, 2, 3, 5, 8, 36, 163, 179, 134, 196, 145)_{7754}$

$(0, 1, 2, 3, 5, 12, 40, 93, 150, 158, 134, 163)_{7840}$

$(0, 1, 2, 3, 5, 12, 73, 155, 133, 24, 115, 161)_{7845}$

$(0, 1, 2, 3, 6, 12, 74, 183, 101, 164, 94, 45)_{7847}$
\ADFvfyParStart{\{217, \{217, \{\{1, 217, 3, 25, \{\{217, 1\}\}\}\}\}\}}

{Order 289}, $\omega = 110$:


$(0, 1, 2, 3, 4, 6, 40, 99, 56, 232, 173, 211)_{201}$

$(0, 1, 2, 3, 4, 6, 88, 184, 48, 225, 129, 27)_{1513}$

$(0, 1, 2, 3, 5, 9, 43, 59, 37, 144, 99, 172)_{3470}$

$(0, 1, 2, 3, 4, 7, 27, 135, 279, 192, 170, 243)_{6383}$

$(0, 1, 2, 3, 4, 8, 67, 11, 88, 245, 182, 18)_{6397}$

$(0, 1, 2, 3, 4, 8, 12, 224, 96, 107, 264, 242)_{6401}$

$(0, 1, 2, 3, 4, 8, 10, 160, 205, 241, 174, 273)_{6406}$

$(0, 1, 2, 3, 4, 8, 80, 260, 187, 245, 169, 69)_{6408}$

$(0, 1, 2, 3, 4, 8, 21, 188, 112, 123, 48, 88)_{6713}$

$(0, 1, 2, 3, 5, 7, 23, 277, 240, 132, 117, 206)_{6753}$

$(0, 1, 2, 3, 5, 8, 24, 240, 136, 123, 133, 90)_{7677}$

$(0, 1, 2, 3, 5, 8, 17, 219, 150, 278, 254, 108)_{7700}$

$(0, 1, 2, 3, 5, 8, 17, 232, 19, 96, 210, 171)_{7754}$

$(0, 1, 2, 3, 5, 10, 34, 59, 127, 19, 228, 212)_{7840}$

$(0, 1, 2, 3, 5, 10, 29, 49, 214, 154, 196, 259)_{7845}$

$(0, 1, 2, 3, 6, 12, 108, 266, 227, 13, 157, 33)_{7847}$
\ADFvfyParStart{\{289, \{289, \{\{1, 289, 4, 110, \{\{289, 1\}\}\}\}\}\}}

{Order 433}, $\omega = 64$:


$(0, 1, 2, 3, 4, 6, 19, 140, 157, 266, 32, 208)_{201}$

$(0, 1, 2, 3, 4, 6, 20, 85, 187, 401, 342, 70)_{1513}$

$(0, 1, 2, 3, 5, 9, 26, 58, 183, 419, 145, 240)_{3470}$

$(0, 1, 2, 3, 4, 7, 14, 155, 207, 317, 393, 147)_{6383}$

$(0, 1, 2, 3, 4, 8, 12, 216, 171, 34, 133, 97)_{6397}$

$(0, 1, 2, 3, 4, 8, 11, 250, 393, 222, 342, 67)_{6401}$

$(0, 1, 2, 3, 4, 8, 10, 221, 290, 380, 262, 240)_{6406}$

$(0, 1, 2, 3, 4, 8, 17, 103, 188, 79, 394, 298)_{6408}$

$(0, 1, 2, 3, 4, 8, 16, 336, 302, 394, 265, 206)_{6713}$

$(0, 1, 2, 3, 5, 7, 15, 169, 63, 195, 33, 393)_{6753}$

$(0, 1, 2, 3, 5, 8, 19, 265, 327, 138, 43, 53)_{7677}$

$(0, 1, 2, 3, 5, 8, 19, 394, 206, 290, 282, 191)_{7700}$

$(0, 1, 2, 3, 5, 8, 16, 44, 152, 171, 334, 222)_{7754}$

$(0, 1, 2, 3, 5, 10, 20, 122, 247, 202, 391, 111)_{7840}$

$(0, 1, 2, 3, 5, 10, 33, 185, 108, 136, 341, 410)_{7845}$

$(0, 1, 2, 3, 6, 12, 14, 339, 300, 69, 160, 253)_{7847}$
\ADFvfyParStart{\{433, \{433, \{\{1, 433, 6, 64, \{\{433, 1\}\}\}\}\}\}}

{Order 577}, $\omega = 27$:


$(0, 1, 2, 3, 4, 6, 14, 501, 69, 300, 402, 539)_{201}$

$(0, 1, 2, 3, 4, 6, 10, 115, 283, 323, 210, 348)_{1513}$

$(0, 1, 2, 3, 5, 9, 15, 315, 190, 509, 115, 250)_{3470}$

$(0, 1, 2, 3, 4, 7, 14, 60, 533, 231, 209, 445)_{6383}$

$(0, 1, 2, 3, 4, 8, 12, 360, 52, 371, 494, 298)_{6397}$

$(0, 1, 2, 3, 4, 8, 14, 348, 419, 295, 34, 212)_{6401}$

$(0, 1, 2, 3, 4, 8, 10, 348, 62, 115, 201, 322)_{6406}$

$(0, 1, 2, 3, 4, 8, 17, 307, 176, 565, 47, 141)_{6408}$

$(0, 1, 2, 3, 4, 8, 16, 360, 308, 210, 283, 53)_{6713}$

$(0, 1, 2, 3, 5, 7, 14, 49, 310, 179, 198, 331)_{6753}$

$(0, 1, 2, 3, 5, 8, 16, 105, 266, 34, 421, 288)_{7677}$

$(0, 1, 2, 3, 5, 8, 16, 181, 256, 50, 139, 37)_{7700}$

$(0, 1, 2, 3, 5, 8, 16, 272, 231, 477, 549, 452)_{7754}$

$(0, 1, 2, 3, 5, 12, 17, 253, 75, 478, 42, 465)_{7840}$

$(0, 1, 2, 3, 5, 12, 17, 417, 356, 479, 158, 249)_{7845}$

$(0, 1, 2, 3, 6, 12, 14, 271, 51, 191, 114, 237)_{7847}$
\ADFvfyParStart{\{577, \{577, \{\{1, 577, 8, 27, \{\{577, 1\}\}\}\}\}\}}

{Order 1009}, $\omega = 139$:


$(0, 1, 2, 3, 4, 6, 13, 982, 338, 658, 314, 547)_{201}$

$(0, 1, 2, 3, 4, 6, 13, 76, 538, 779, 663, 978)_{1513}$

$(0, 1, 2, 3, 5, 10, 14, 277, 428, 808, 934, 949)_{3470}$

$(0, 1, 2, 3, 4, 7, 14, 79, 88, 182, 965, 738)_{6383}$

$(0, 1, 2, 3, 4, 8, 14, 100, 62, 530, 986, 322)_{6397}$

$(0, 1, 2, 3, 4, 8, 14, 45, 771, 428, 924, 652)_{6401}$

$(0, 1, 2, 3, 4, 8, 17, 55, 80, 117, 659, 232)_{6406}$

$(0, 1, 2, 3, 4, 8, 17, 54, 947, 123, 468, 85)_{6408}$

$(0, 1, 2, 3, 4, 8, 17, 182, 463, 312, 866, 597)_{6713}$

$(0, 1, 2, 3, 5, 7, 14, 59, 240, 979, 613, 86)_{6753}$

$(0, 1, 2, 3, 5, 8, 17, 115, 330, 154, 22, 488)_{7677}$

$(0, 1, 2, 3, 5, 8, 17, 68, 121, 481, 425, 189)_{7700}$

$(0, 1, 2, 3, 5, 8, 17, 55, 987, 894, 658, 855)_{7754}$

$(0, 1, 2, 3, 5, 14, 19, 33, 320, 363, 166, 523)_{7840}$

$(0, 1, 2, 3, 5, 14, 19, 50, 351, 842, 86, 40)_{7845}$

$(0, 1, 2, 3, 6, 12, 14, 33, 132, 231, 524, 972)_{7847}$
\ADFvfyParStart{\{1009, \{1009, \{\{1, 1009, 14, 139, \{\{1009, 1\}\}\}\}\}\}}
%
\end{proof}


\section{Graph designs: general constructions}\label{sec:Graph-designs-general-constructions}
In Propositions~\ref{prp:24-4}--\ref{prp:18^5-9^9} we describe some general constructions for 12-vertex, 6-regular graph designs.
We refer to Lemma~\ref{lem-GDD-existence} for the existence of the various group divisible designs mentioned.
Observe that a $G$-design of order 1 always exists---it is the empty set.
In what follows we tacitly assume this trivial case.

%
\begin{proposition}\label{prp:24-4}
Let $G$ be a 12-vertex, 6-regular graph.
Suppose there exists a $G$-decomposition of the complete multipartite graph $K_{24^4}$.
Suppose also that there exist $G$-designs of orders $73$, $145$, $217$ and $433$.
Then there exist $G$-designs of order $n$ for all positive integers $n \equiv 1 \adfmod{72}$.
\end{proposition}
%
\begin{proof}
Let $t \ge 5$ and $u \ge 4$ be integers.

Take a 4-GDD of type $6^t$,
inflate its points by a factor of 24 and
replace its blocks by $G$-decompositions of $K_{24^4}$.
Add a new point and overlay each group plus the new point with a $G$-design of order 145.
The result is a $G$-design of order $144t + 1$ for $t \ge 5$.

Take a 4-GDD of type $6^u 3^1$,
inflate its points by a factor of 24 and
replace its blocks by $G$-decompositions of $K_{24^4}$.
Add a new point and overlay each group plus the new point with a $G$-design of order 73 or 145, as appropriate.
The result is a $G$-design of order $144u + 73$ for $u \ge 4$.

We deal with the orders missed, namely 289, 361, 505 and 577, by similar constructions.
For brevity we just indicate the ingredients.

For order 289,
use a 4-GDD of type $3^4$ with $G$-decompositions of $K_{24^4}$ and $K_{73}$.

For order 361,
use a 4-GDD of type $3^5$ with $G$-decompositions of $K_{24^4}$ and $K_{73}$.

For order 505,
use a 4-GDD of type $3^5 6^1$ with $G$-decompositions of $K_{24^4}$, $K_{73}$ and $K_{145}$.

For order 577,
use a 4-GDD of type $3^8$ with $G$-decompositions of $K_{24^4}$ and $K_{73}$.
\end{proof}
%
\begin{proposition}\label{prp:18^5-18^4}
Let $G$ be a 12-vertex, 6-regular graph.
Suppose there exist $G$-decompositions of $K_{18^4}$ and $K_{18^5}$.
Suppose also that there exist $G$-designs of orders $73$, $145$, $217$ and $433$.
Then there exist $G$-designs of order $n$ for all positive integers $n \equiv 1 \adfmod{72}$.
\end{proposition}
%
\begin{proof}
Let $m \in \{0,4,8\}$ and $t \ge \max\{1, m/4\}$ be integers.
Take a $\{4,5\}$-GDD of type $4^{3t + 1} m^1$,
inflate its points by a factor of 18 and
replace its blocks by $G$-decompositions of $K_{18^4}$ or $K_{18^5}$, as appropriate.
Add a new point and overlay each group plus the new point with a $G$-design of order 73 or 145, as appropriate.
The result is
a $G$-design of order $216t + 18m + 73$ for $m \in \{0,4,8\}$, $t \ge \max\{1, m/4\}$.
No further constructions are needed.
%
\end{proof}
%
\begin{proposition}\label{prp:18^5-9^9}
Let $G$ be a 12-vertex, 6-regular graph.
Suppose there exist $G$-decompositions of $K_{18^5}$, $K_{9^9}$ and either $K_{6^7}$ or $K_{72^7}$.
Suppose also that there exist $G$-designs of orders $73$, $145$, $217$, $289$, $577$ and $1009$.
Then there exist $G$-designs of order $n$ for all positive integers $n \equiv 1 \adfmod{72}$.
\end{proposition}
%
\begin{proof}
Let $m \in \{0, 4, 8, 12, 16\}$ and $t$ be integers such that
\begin{equation}
\label{eqn:18^5-9^9}
t \ge \left\{\begin{array}{l} 1 ~ \text{if } m \in \{0, 4\}, \\
                              3 ~ \text{if } m = 8, \\
                              2 ~ \text{if } m = 12, \\
                              4 ~ \text{if } m = 16.
               \end{array} \right.
\end{equation}
Take a 5-GDD of type $4^{5t} m^1$,
inflate its points by a factor of 18 and
replace its blocks by $G$-decompositions of $K_{18^5}$.
Add a new point and overlay each group plus the new point with a $G$-design of order 73, 145, 217 or 289, as appropriate.
The result is
a $G$-design of order $360t + 18m + 1$ for $m \in \{0, 4, 8, 12, 16\}$ and $t$ satisfying (\ref{eqn:18^5-9^9}).

We deal with the orders missed, namely 505, 649, 865, 1369, by similar constructions.
If the multipartite graph is $K_{f^e}$, we inflate each point of the GDD by a factor of $f$.

For order 505, either
use a 7-GDD of type $12^7$ with $G$-decompositions of $K_{6^7}$ and $K_{73}$, or
use a 7-GDD of type $1^7$ with $G$-decompositions of $K_{72^7}$ and $K_{73}$.

For order 649,
use a 9-GDD of type $8^9$ with $G$-decompositions of $K_{9^9}$ and $K_{73}$.

For order 865,
use a 5-GDD of type $8^6$ with $G$-decompositions of $K_{18^5}$ and $K_{145}$.

For order 1369,
use a 5-GDD of type $12^5 16^1$ with $G$-decompositions of $K_{18^5}$, $K_{217}$ and $K_{289}$.
\end{proof}

Now we are ready to prove our main theorems.
%
\begin{theorem}\label{thm:spectra-complete-main}
For graphs
$201$, $1513$, $3470$, $6383$, $6397$, $6401$, $6406$, $6408$, $6713$, $6753$, $7677$, $7700$, $7754$, $7840$, $7845$ and $7847$,
a design of order $n$ exists if and only if $n \equiv 1 \adfmod{72}$.
\end{theorem}
%
\begin{proof}
For
201, 6383, 6397, 6401, 6753, 7677, 7754,
use Proposition \ref{prp:24-4} with a decomposition of
$K_{24^4}$
and design orders
73, 145, 217, 433.

For
6406, 6408, 7845,
use Proposition \ref{prp:18^5-9^9} with decompositions of
$K_{18^5}$, $K_{6^7}$, $K_{9^9}$
and design orders
73, 145, 217, 289, 577, 1009.

For
1513, 3470, 6713,
use Proposition \ref{prp:18^5-18^4} with decompositions of
$K_{18^4}$, $K_{18^5}$
and design orders
73, 145, 217, 433.

For
7700, 7840, 7847,
use Proposition \ref{prp:18^5-9^9} with decompositions of
$K_{18^5}$, $K_{72^7}$, $K_{9^9}$
and design orders
73, 145, 217, 289, 577, 1009.

See Lemmas~\ref{lem-Design-24-4}--\ref{lem-Design-73-et-al} for the relevant graph decompositions.
\end{proof}

As explained in the Introduction, Theorem~\ref{thm:spectra-complete-main} deals only with the graphs for which we
have provided decomposition details in Lemmas~\ref{lem-Design-24-4}--\ref{lem-Design-73-et-al}.
The next theorem represents all of our successful design spectrum completions.
%
\begin{theorem}\label{thm:spectra-complete-appendix}
%
%
For $7788$ 12-vertex, 6-regular graphs, including the 5-colourable vertex-transitive graphs and all of the 3-chromatic graphs,
$1$, $2$, $3$, $4$, $17$, $18$, $20$, $22$, $23$, $24$, $201$, $203$, $206$, $207$, $228$, $312$, $527$, $529$, $590$, $599$, $601$, $850$, $1106$, $1233$, $1261$, $1698$, $1702$, $1825$, $1835$, $1839$, $2040$, $2045$, $2051$, $2053$, $2471$, $2562$, $2563$, $2574$, $2581$, $3179$, $3191$, $3193$, $3241$, $3243$, $6383$, $6385$, $6390$, $6397$ and $6401$,
a design of order $n$ exists if and only if $n \equiv 1 \adfmod{72}$.
\end{theorem}
%
\begin{proof}
(i) In Part A of the appendix we give decomposition details of
$$ K_{24^4}, K_{73}, K_{145}, K_{217}, K_{433}$$
for 6311 graphs, including seven covered by Theorem~\ref{thm:spectra-complete-main}.
Use Proposition~\ref{prp:24-4}.

(ii) In Part B of the appendix we give decomposition details of
$$ K_{18^5}, K_{6^7}, K_{9^9}, K_{73}, K_{145}, K_{217}, K_{289}, K_{577}, K_{1009}$$
for 1471 graphs, including three covered by Theorem~\ref{thm:spectra-complete-main}.
Use Proposition~\ref{prp:18^5-9^9}.

The six graphs, 1513, 3470, 6713, 7700, 7840, 7847, not covered by (i) and (ii) are dealt with in Theorem~\ref{thm:spectra-complete-main}.
\end{proof}
%
\begin{theorem}\label{thm:spectra-missing-505-appendix}
For eleven 4-chromatic graphs,
$10$, $13$, $59$, $130$, $211$, $432$, $551$, $3281$, $6729$, $7679$, $7743$,
and forty-three 5-chromatic graphs,
$16$, $163$, $424$, $635$, $659$, $670$, $671$, $687$, $692$, $701$, $702$, $707$, $722$, $733$, $1063$, $1438$, $3101$, $3443$, $3447$, $4001$, $4069$, $4070$, $4074$, $4096$, $4108$, $4317$, $4764$, $4778$, $5701$, $5859$, $5913$, $6339$, $6391$, $6657$, $6751$, $7353$, $7421$, $7531$, $7603$, $7667$, $7752$, $7761$, $7803$,
a design of order $n$ exists if and only if $n \equiv 1 \adfmod{72}$,
with the possible exception of $n = 505$.
\end{theorem}
%
\begin{proof}
For each of the 54 graphs in the statement of the theorem, we have decompositions of
$$ K_{18^5}, K_{9^9}, K_{73}, K_{145}, K_{217}, K_{289}, K_{577}, K_{1009}.$$
The details are in Part C of the appendix.

Use Proposition \ref{prp:18^5-9^9} omitting the construction of a design of order 505
since we do not have decompositions of $K_{6^7}$ or $K_{72^7}$ for any of the graphs.
\end{proof}


\section{Concluding remarks}\label{sec:Concluding-remarks}

We cannot help thinking that fate has been kind to us.
Prior to carrying out the work for this paper,
we would have been forgiven for believing that any attempt to obtain
the design spectrum of a 6-regular graph with chromatic number greater than 2 would be doomed.
Even for the smallest example, $K_7$, there are unresolved cases---twenty-one stated in \cite[Table 1]{AdamsBryantBuchanan2008}.
And yet we have in our paper completely solved the spectrum problem for thousands of 6-regular graphs, a substantial number of them 5-chromatic.
Here we offer some explanations.

The necessary condition for a design of order $n$ for 12 vertices and 6-regularity is particularly simple, $n \equiv 1 \adfmod{72}$.
It is well known to graph-design theorists that residue class 1 modulo $2|E(G)|$ is by far the easiest.
Take the truncated cuboctahedron, for example. In 2013 Forbes \& Griggs published the
solution only for design orders $n \equiv 1 \adfmod{144}$, \cite{ForbesGriggs2013}.
However, $n$ can belong to another residue class, $64 \adfmod{144}$, and
it took four extra years plus one extra author (T.\ J.\ Forbes) to resolve this case, \cite{ForbesForbesGriggs2017}.

Other things being equal, 5-chromatic graphs are usually much harder to process than 4-colourable ones.
Constructions like Proposition~\ref{prp:18^5-9^9} require a suitable infinite supply of 5-GDDs.
In contrast to 2-, 3- and 4-GDDs, group divisible designs with block size 5 are scarce,
and in most cases the right ones for a particular design spectrum are unavailable.
However, the publication of the paper by Abel, Ge, Greig \& Ling, \cite{AbelGeGreigLing2009}, is most fortunate.
They obtain substantial results concerning the existence of $(v, \{5, w^*\}, 1)$-PBDs and hence of 5-GDDs of type $4^u m^1$,
which turn out to be precisely what we want to combine with $G$-decompositions of $K_{18^5}$.

Obtaining a direct construction of a 6-regular graph decomposition by backtracking combined with random processes
is in general essentially hopeless unless there is an automorphism of large order.
We were therefore pleasantly surprised when we discovered that all of our $G$-decompositions except $K_{72^7}$
can be obtained from single base blocks.

On the other hand, there was a limit to our good fortune.
Recall from Theorem~\ref{thm:spectra-missing-505-appendix} that we were obliged to accept the possible exception of order 505
for a specific set $\mathcal{C}$, say, of 54 graphs.
There are various ways to try to obtain this design order.
\begin{center}
\begin{tabular}{c|c|c}
GDD & decompositions & design orders \\
\hline
4-GDD type $3^5 6^1$   & $K_{24^4}$             & \rule{0mm}{4mm} 73, 145 \rule{0mm}{4mm}\\
4-GDD type $4^7$       & $K_{18^4}$             & 73 \\
7-GDD type $1^7$       & $K_{72^7}$             & 73 \\
7-GDD type $12^7$      & $K_{6^7}$              & 73
\end{tabular}
\end{center}
The first two options won't work for 5-chromatic graphs,
the $K_{72^7}$ option did not achieve very much, and
we suspect a $G$-decomposition of $K_{6^7}$ generated from a single base block does not exist for any $G \in \mathcal{C}$.
However, we had no difficulty obtaining $G$-decompositions of $K_{12^7}$ for all $G \in \mathcal{C}$.
With these in place there appears to be an obvious way to construct a $G$-design of order 505.
Use the decomposition of $K_{12^7}$ with design order 73 and
a 7-GDD of type $6^7$ created by removing a point from a projective plane of order 6.


\section*{Acknowledgement}

We would like to thank Dr Markus Meringer for creating and making available
the edge sets of the 7848 connected 12-vertex, 5-regular graphs, \cite{Meringer1999,MeringerReggraphs1999}.

\section*{ORCID}

\noindent A. D. Forbes     \url{https://orcid.org/0000-0003-3805-7056} \\
C. G. Rutherford           \url{https://orcid.org/0000-0003-1924-207X}


\adfhide
{
\section*{Declarations}

\subsection*{Conflicts}

There are no conflicts of interest.

\subsection*{Funding}

No funds, grants, or other support were received during the preparation of this manuscript.

\subsection*{Data availability}

We do not analyse or generate any datasets.
}


\newcommand{\adfgw}{\textwidth}
\begin{figure}
\caption{Vertex-transitive graphs}
\begin{center}
\includegraphics[width=\adfgw,trim=0mm 0mm 0mm 0mm, clip]{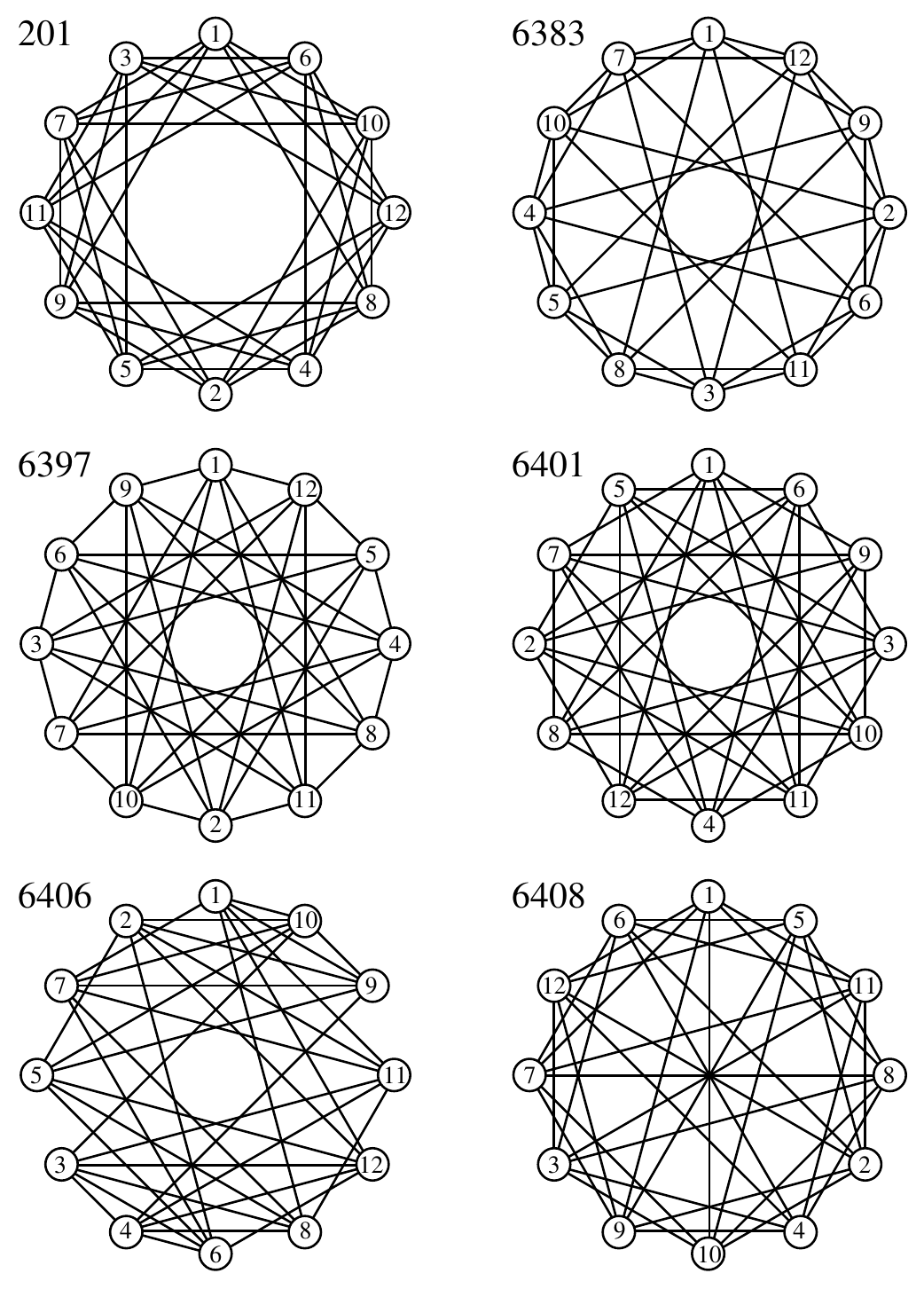} 
\end{center}
\label{fig:graph-pictures-first}
\end{figure}

\begin{figure}
\caption{Vertex-transitive graphs}
\begin{center}
\includegraphics[width=\adfgw,trim=0mm 0mm 0mm 0mm, clip]{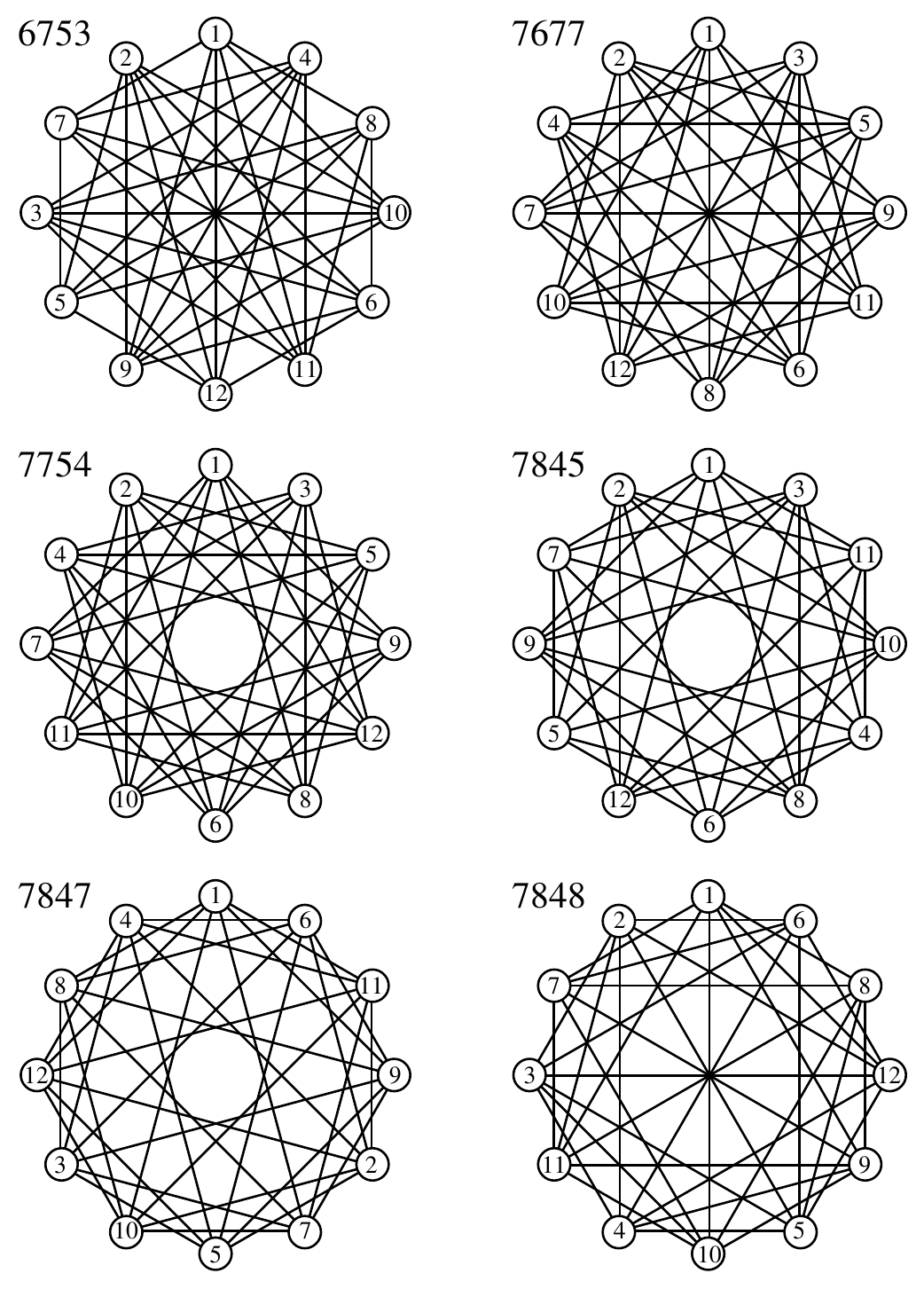} 
\end{center}
\end{figure}

\begin{figure}
\caption{Graphs}
\begin{center}
\includegraphics[width=\adfgw,trim=0mm 0mm 0mm 0mm, clip]{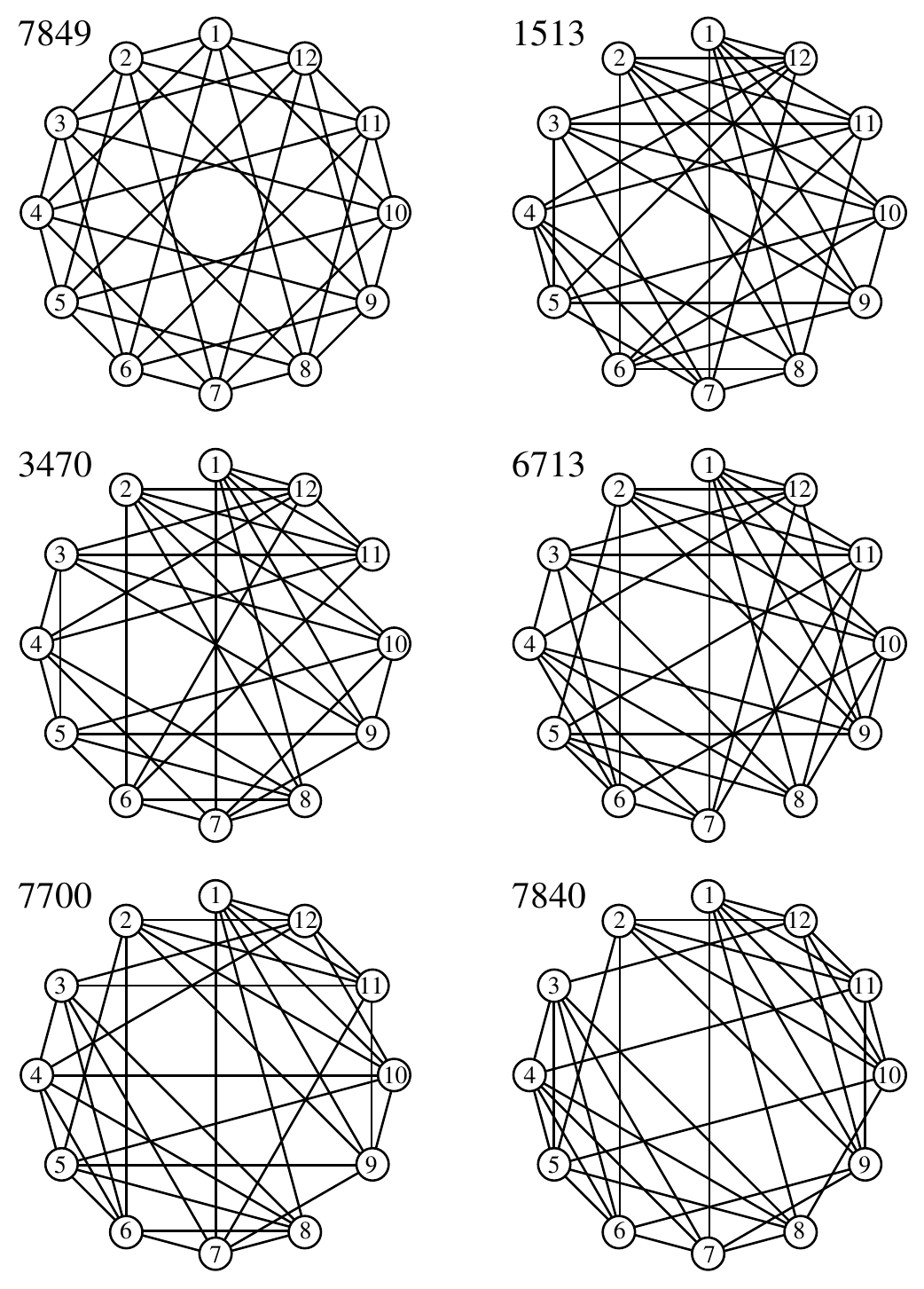} 
\end{center}
\label{fig:graph-pictures-last}
\end{figure}


\newpage
\hoffset -15mm \voffset -15mm
\topmargin 0mm
\headheight 5mm
\headsep 5mm
\textheight 250mm
\footskip 0mm
\evensidemargin 10mm
\oddsidemargin  10mm
\textwidth 180mm

\appendix

\noindent {\bf APPENDIX}\\~

\parindent 0mm

\newcommand{\adfSep}{\hspace{4mm}}

\noindent The details for a 12-vertex, 6-regular graph are coded as
$$\text{\boldmath $g$}\adfSep (a_1, a_2, \dots, a_{11})\adfSep d_1\adfSep B_1\adfSep d_2\adfSep B_2\adfSep \dots,$$
where:

\vskip 2mm
\begin{enumerate}
\setlength{\itemsep}{2mm}
\item[$g$] is the graph number;
\item[$a_i$] is row $i$ of the above-diagonal part of the graph $g$'s adjacency matrix,
    treated as a binary number and converted to decimal;
\item[$d_j$] specifies the design type represented by $B_j$; either \\[2mm]
    $d_j$ is a pair $(p_j,u_j)$ to indicate a decomposition into copies of graph $g$ of
    the complete multipartite graph with $u_j$ parts of size $p_j$, or \\[2mm]
    $d_j$ is a single number to indicate a (graph $g$)-design of order $d_j$;
\item[$B_j$] is the single base block for the design of type $d_j$.
\end{enumerate}

\vskip 3mm
The parameters $(a_1, a_2, \dots, a_{11})$ for graphs 1--7848 are derived from
the complements of 12-vertex, 5-regular graphs in Meringer's list, \cite{MeringerReggraphs1999}.

\vskip 3mm
The construction of the decomposition from the base block $B_j$ depends only on the design type $d_j$, as explained in
Lemmas~\ref{lem-Design-24-4},
\ref{lem-Design-18-5},
\ref{lem-Design-6-7},
\ref{lem-Design-9-9} and \ref{lem-Design-73-et-al}.

\vskip 3mm
There are three parts:

\vskip 2mm
\begin{enumerate}
\setlength{\itemsep}{2mm}
\item[A:] decomposition of $K_{24^4}$ and designs of orders 73, 145, 217, 433 for 6311 graphs;
          page~\pageref{app:PART-A};
\item[B:] decompositions of $K_{18^5}$, $K_{6^7}$, $K_{9^9}$ and designs of orders 73, 145, 217, 289, 577, 1009 for 1471 graphs;
          page~\pageref{app:PART-B};
\item[C:] decompositions of $K_{18^5}$, $K_{9^9}$ and designs of orders 73, 145, 217, 289, 577, 1009 for 54 graphs;
          page~\pageref{app:PART-C}.
\end{enumerate}

$~$\newpage

\tiny \parskip 0.5mm \parindent 0mm


{\normalsize\bf PART A}\\~\label{app:PART-A}

\begin{tabular}{@{}llllll}
{\bf 1} & (63,63,63,63,95,47,16,0,0,0,0) &
(24,4) & (75,43,19,83,54,59,8,1,13,44,56,92)  &
73 & (0,1,2,3,4,5,59,41,7,13,25,31) \\
145 & (0,1,2,3,4,5,7,83,18,30,41,54)  &
217 & (0,1,2,3,4,5,54,37,63,83,89,108)  &
433 & (0,1,2,3,4,6,68,335,188,240,299,430)
\end{tabular}

\begin{tabular}{@{}llllll}
{\bf 2} & (63,63,63,95,95,39,24,0,0,0,0) &
(24,4) & (84,40,8,36,72,30,3,10,26,5,85,93)  &
73 & (0,1,2,3,4,10,41,5,60,48,53,65) \\
145 & (0,1,2,3,4,5,78,18,65,30,41,53)  &
217 & (0,1,2,3,4,5,189,23,173,8,73,98)  &
433 & (0,1,2,3,4,5,67,50,299,33,150,197)
\end{tabular}

\begin{tabular}{@{}llllll}
{\bf 3} & (63,63,63,95,111,23,24,0,0,0,0) &
(24,4) & (42,35,43,49,50,51,81,0,72,32,56,76)  &
73 & (0,1,2,3,4,5,36,45,57,8,15,27) \\
145 & (0,1,2,3,4,5,10,71,23,36,48,60)  &
217 & (0,1,2,3,4,5,12,114,63,39,81,89)  &
433 & (0,1,2,3,4,5,124,287,50,85,160,197)
\end{tabular}

\begin{tabular}{@{}llllll}
{\bf 4} & (63,63,63,95,111,51,4,8,0,0,0) &
(24,4) & (47,87,67,38,58,73,68,76,21,81,16,40)  &
73 & (0,1,2,3,4,5,8,59,41,13,25,31) \\
145 & (0,1,2,3,4,5,10,53,107,119,18,66)  &
217 & (0,1,2,3,4,5,8,97,34,198,49,74)  &
433 & (0,1,2,3,4,5,9,119,329,240,44,188)
\end{tabular}

\begin{tabular}{@{}llllll}
{\bf 5} & (63,63,63,95,111,51,16,0,4,0,0) &
(24,4) & (44,76,52,80,48,69,35,74,17,90,51,95)  &
73 & (0,1,2,3,4,5,53,65,7,30,19,36) \\
145 & (0,1,2,3,4,5,72,65,17,120,41,53)  &
217 & (0,1,2,3,4,5,14,141,38,89,63,107)  &
433 & (0,1,2,3,4,5,85,197,124,287,49,160)
\end{tabular}

\begin{tabular}{@{}llllll}
{\bf 6} & (63,63,63,207,79,51,20,8,0,0,0) &
(24,4) & (77,29,53,70,69,95,58,60,47,15,20,80)  &
73 & (0,1,2,3,4,6,20,58,8,64,28,43) \\
145 & (0,1,2,3,4,7,75,112,96,88,83,124)  &
217 & (0,1,2,3,4,6,19,29,101,150,48,95)  &
433 & (0,1,2,3,4,6,14,266,317,408,46,208)
\end{tabular}

\begin{tabular}{@{}llllll}
{\bf 7} & (63,63,63,207,87,27,28,0,0,0,0) &
(24,4) & (34,37,38,24,78,85,28,59,51,75,39,67)  &
73 & (0,1,2,3,4,6,44,28,69,36,15,56) \\
145 & (0,1,2,3,4,6,15,108,120,87,71,98)  &
217 & (0,1,2,3,4,6,15,212,34,114,57,138)  &
433 & (0,1,2,3,4,6,11,206,406,388,302,338)
\end{tabular}

\begin{tabular}{@{}llllll}
{\bf 9} & (63,63,63,207,87,43,20,8,0,0,0) &
(24,4) & (24,59,92,8,35,2,41,22,1,70,49,89)  &
73 & (0,1,2,3,4,11,28,59,40,71,52,64) \\
145 & (0,1,2,3,4,9,54,71,38,65,50,122)  &
217 & (0,1,2,3,4,6,24,91,72,32,123,175)  &
433 & (0,1,2,3,4,6,24,303,363,298,126,427)
\end{tabular}

\begin{tabular}{@{}llllll}
{\bf 11} & (63,63,63,207,87,57,6,8,0,0,0) &
(24,4) & (43,72,8,82,36,23,50,38,61,17,33,37)  &
73 & (0,1,2,3,4,11,5,25,47,20,35,13) \\
145 & (0,1,2,3,4,6,29,129,54,42,139,66)  &
217 & (0,1,2,3,4,6,29,102,94,45,84,152)  &
433 & (0,1,2,3,4,6,10,277,388,103,352,126)
\end{tabular}

\begin{tabular}{@{}llllll}
{\bf 12} & (63,63,63,207,87,57,10,4,0,0,0) &
(24,4) & (10,30,22,38,29,52,15,87,9,40,84,55)  &
73 & (0,1,2,3,4,6,46,36,15,25,20,69) \\
145 & (0,1,2,3,4,6,58,81,48,22,35,140)  &
217 & (0,1,2,3,4,6,29,19,118,156,87,39)  &
433 & (0,1,2,3,4,6,10,126,333,160,352,388)
\end{tabular}

\begin{tabular}{@{}llllll}
{\bf 14} & (63,63,63,207,115,60,1,2,4,0,0) &
(24,4) & (54,38,94,46,75,92,33,81,61,95,56,24)  &
73 & (0,1,2,3,4,11,32,40,5,25,13,58) \\
145 & (0,1,2,3,4,6,61,78,16,108,97,30)  &
217 & (0,1,2,3,4,6,55,118,62,133,183,208)  &
433 & (0,1,2,3,4,6,10,417,93,126,221,253)
\end{tabular}

\begin{tabular}{@{}llllll}
{\bf 15} & (63,63,63,207,115,60,1,4,2,0,0) &
(24,4) & (81,1,89,65,76,10,75,63,20,60,74,18)  &
73 & (0,1,2,3,4,10,41,30,71,65,25,60) \\
145 & (0,1,2,3,4,6,52,77,40,16,66,60)  &
217 & (0,1,2,3,4,6,34,206,104,153,128,41)  &
433 & (0,1,2,3,4,6,10,93,417,126,221,253)
\end{tabular}

\begin{tabular}{@{}llllll}
{\bf 17} & (63,63,95,95,111,19,28,0,0,0,0) &
(24,4) & (14,41,52,80,49,58,12,27,11,59,75,91)  &
73 & (0,1,2,3,4,10,28,53,59,72,36,64) \\
145 & (0,1,2,3,4,5,131,78,41,120,52,64)  &
217 & (0,1,2,3,4,5,32,212,48,198,74,100)  &
433 & (0,1,2,3,4,5,127,175,9,197,85,157)
\end{tabular}

\begin{tabular}{@{}llllll}
{\bf 18} & (63,63,95,95,111,35,20,8,0,0,0) &
(24,4) & (62,78,6,47,22,16,7,28,85,76,17,65)  &
73 & (0,1,2,3,4,5,22,26,46,41,10,16) \\
145 & (0,1,2,3,4,5,72,11,126,114,48,60)  &
217 & (0,1,2,3,4,5,17,68,45,94,74,102)  &
433 & (0,1,2,3,4,5,66,240,278,386,85,151)
\end{tabular}

\begin{tabular}{@{}llllll}
{\bf 19} & (63,63,95,95,111,49,18,0,4,0,0) &
(24,4) & (57,61,29,5,81,54,64,95,66,80,18,83)  &
73 & (0,1,2,3,4,5,38,47,20,52,8,15) \\
145 & (0,1,2,3,4,5,139,54,18,66,119,42)  &
217 & (0,1,2,3,4,5,32,212,22,122,173,74)  &
433 & (0,1,2,3,4,5,170,367,138,248,8,102)
\end{tabular}

\begin{tabular}{@{}llllll}
{\bf 20} & (63,63,95,111,119,11,20,8,0,0,0) &
(24,4) & (71,51,87,67,12,89,42,1,50,9,6,94)  &
73 & (0,1,2,3,4,10,63,51,30,72,41,58) \\
145 & (0,1,2,3,4,5,68,131,120,62,42,53)  &
217 & (0,1,2,3,4,5,20,145,96,71,47,102)  &
433 & (0,1,2,3,4,5,98,229,171,27,33,140)
\end{tabular}

\begin{tabular}{@{}llllll}
{\bf 21} & (63,63,95,111,119,25,6,8,0,0,0) &
(24,4) & (32,44,4,24,84,69,93,74,7,11,43,26)  &
73 & (0,1,2,3,4,5,14,20,41,26,31,9) \\
145 & (0,1,2,3,4,5,52,65,78,29,108,20)  &
217 & (0,1,2,3,4,6,21,195,130,123,177,149)  &
433 & (0,1,2,3,4,5,8,157,171,208,413,101)
\end{tabular}

\begin{tabular}{@{}llllll}
{\bf 22} & (63,63,95,111,119,25,10,4,0,0,0) &
(24,4) & (32,36,4,12,40,9,14,75,59,94,3,19)  &
73 & (0,1,2,3,4,5,8,20,14,30,25,36) \\
145 & (0,1,2,3,4,5,66,18,141,53,120,43)  &
217 & (0,1,2,3,4,5,18,145,127,102,69,42)  &
433 & (0,1,2,3,4,5,11,229,278,265,30,101)
\end{tabular}

\begin{tabular}{@{}llllll}
{\bf 23} & (63,63,95,111,119,25,24,0,0,2,0) &
(24,4) & (49,57,25,1,29,64,59,30,46,66,72,71)  &
73 & (0,1,2,3,4,10,47,25,19,5,36,63) \\
145 & (0,1,2,3,4,5,66,43,18,9,54,30)  &
217 & (0,1,2,3,4,5,87,81,114,34,42,59)  &
433 & (0,1,2,3,4,5,12,118,190,73,111,47)
\end{tabular}

\begin{tabular}{@{}llllll}
{\bf 24} & (63,63,95,111,119,56,1,2,4,0,0) &
(24,4) & (85,41,21,13,65,94,68,40,3,78,74,19)  &
73 & (0,1,2,3,4,5,14,36,20,8,52,47) \\
145 & (0,1,2,3,4,5,9,70,89,114,47,22)  &
217 & (0,1,2,3,4,5,10,198,188,48,98,148)  &
433 & (0,1,2,3,4,5,12,118,328,44,73,190)
\end{tabular}

\begin{tabular}{@{}llllll}
{\bf 25} & (63,63,95,111,119,56,1,8,0,2,0) &
(24,4) & (45,77,81,13,49,16,90,54,75,91,92,20)  &
73 & (0,1,2,3,4,5,69,15,41,25,31,57) \\
145 & (0,1,2,3,4,5,10,96,83,18,119,107)  &
217 & (0,1,2,3,4,5,120,204,128,31,153,178)  &
433 & (0,1,2,3,4,5,52,206,119,12,280,31)
\end{tabular}

\begin{tabular}{@{}llllll}
{\bf 26} & (63,63,95,159,103,43,20,8,0,0,0) &
(24,4) & (46,66,78,67,14,8,77,0,5,68,29,53)  &
73 & (0,1,2,3,4,19,55,46,26,63,35,71) \\
145 & (0,1,2,3,4,6,14,119,72,95,86,110)  &
217 & (0,1,2,3,4,6,14,183,132,83,77,158)  &
433 & (0,1,2,3,4,6,30,227,373,134,102,267)
\end{tabular}

\begin{tabular}{@{}llllll}
{\bf 27} & (63,63,95,159,103,57,6,8,0,0,0) &
(24,4) & (49,85,62,25,51,17,26,52,63,8,56,64)  &
73 & (0,1,2,3,4,33,54,5,27,11,36,16) \\
145 & (0,1,2,3,4,6,95,11,84,30,131,107)  &
217 & (0,1,2,3,4,6,40,14,44,67,94,144)  &
433 & (0,1,2,3,4,6,34,47,135,101,267,208)
\end{tabular}

\begin{tabular}{@{}llllll}
{\bf 28} & (63,63,95,159,103,57,10,4,0,0,0) &
(24,4) & (48,80,72,60,37,26,95,39,29,74,38,55)  &
73 & (0,1,2,3,4,15,65,5,59,35,28,22) \\
145 & (0,1,2,3,4,6,72,16,108,96,30,63)  &
217 & (0,1,2,3,4,6,78,121,95,174,71,199)  &
433 & (0,1,2,3,4,6,13,247,208,43,120,297)
\end{tabular}

\begin{tabular}{@{}llllll}
{\bf 29} & (63,63,95,175,55,59,4,8,0,0,0) &
(24,4) & (45,61,36,20,21,69,42,10,59,27,74,86)  &
73 & (0,1,2,3,4,10,63,19,26,41,31,71) \\
145 & (0,1,2,3,4,5,69,35,11,101,23,61)  &
217 & (0,1,2,3,4,5,17,168,178,120,69,94)  &
433 & (0,1,2,3,4,5,11,70,191,84,131,156)
\end{tabular}

\begin{tabular}{@{}llllll}
{\bf 30} & (63,63,95,175,55,59,8,4,0,0,0) &
(24,4) & (7,75,87,91,12,20,21,25,42,18,10,50)  &
73 & (0,1,2,3,4,10,6,24,58,31,36,64) \\
145 & (0,1,2,3,4,5,9,55,84,107,18,29)  &
217 & (0,1,2,3,4,5,63,138,54,132,38,108)  &
433 & (0,1,2,3,4,5,19,384,215,307,70,185)
\end{tabular}

\begin{tabular}{@{}llllll}
{\bf 31} & (63,63,95,175,87,43,20,8,0,0,0) &
(24,4) & (61,1,38,53,46,9,86,43,88,87,23,40)  &
73 & (0,1,2,3,4,6,19,7,65,36,31,53) \\
145 & (0,1,2,3,4,6,65,131,140,95,42,120)  &
217 & (0,1,2,3,4,6,48,142,98,200,33,73)  &
433 & (0,1,2,3,4,6,11,82,122,195,280,388)
\end{tabular}

\begin{tabular}{@{}llllll}
{\bf 32} & (63,63,95,175,87,43,24,4,0,0,0) &
(24,4) & (7,87,35,58,83,94,53,14,20,5,52,84)  &
73 & (0,1,2,3,4,7,53,49,9,14,43,58) \\
145 & (0,1,2,3,4,6,11,47,72,91,23,114)  &
217 & (0,1,2,3,4,6,33,16,127,179,144,154)  &
433 & (0,1,2,3,4,6,11,278,109,63,242,263)
\end{tabular}

\begin{tabular}{@{}llllll}
{\bf 33} & (63,63,95,175,87,51,12,8,0,0,0) &
(24,4) & (39,75,42,79,34,19,80,12,89,21,8,36)  &
73 & (0,1,2,3,4,10,12,70,22,35,47,59) \\
145 & (0,1,2,3,4,6,72,55,102,61,37,126)  &
217 & (0,1,2,3,4,6,172,186,98,80,21,73)  &
433 & (0,1,2,3,4,6,9,197,48,85,126,160)
\end{tabular}

\begin{tabular}{@{}llllll}
{\bf 34} & (63,63,95,175,87,51,24,0,4,0,0) &
(24,4) & (73,69,70,85,66,65,4,27,87,36,44,76)  &
73 & (0,1,2,3,4,25,40,18,27,65,6,45) \\
145 & (0,1,2,3,4,6,17,77,84,120,41,96)  &
217 & (0,1,2,3,4,6,163,24,8,99,149,174)  &
433 & (0,1,2,3,4,6,31,135,24,267,208,337)
\end{tabular}

\begin{tabular}{@{}llllll}
{\bf 35} & (63,63,95,175,87,57,6,8,0,0,0) &
(24,4) & (46,25,5,41,51,34,67,0,3,52,60,64)  &
73 & (0,1,2,3,4,11,5,32,40,20,27,14) \\
145 & (0,1,2,3,4,6,123,51,72,64,112,58)  &
217 & (0,1,2,3,4,6,15,127,79,33,42,154)  &
433 & (0,1,2,3,4,6,11,274,69,184,419,384)
\end{tabular}

\begin{tabular}{@{}llllll}
{\bf 36} & (63,63,95,175,87,57,10,4,0,0,0) &
(24,4) & (43,51,86,83,26,75,30,57,60,76,85,48)  &
73 & (0,1,2,3,4,10,48,70,12,43,19,24) \\
145 & (0,1,2,3,4,6,11,59,72,23,102,114)  &
217 & (0,1,2,3,4,6,18,142,121,128,69,43)  &
433 & (0,1,2,3,4,6,11,158,83,243,49,123)
\end{tabular}

\begin{tabular}{@{}llllll}
{\bf 37} & (63,63,95,175,87,57,18,0,4,0,0) &
(24,4) & (14,50,71,58,19,52,87,65,41,16,40,57)  &
73 & (0,1,2,3,4,10,12,70,20,36,26,31) \\
145 & (0,1,2,3,4,6,17,77,84,29,53,41)  &
217 & (0,1,2,3,4,6,11,133,108,64,38,144)  &
433 & (0,1,2,3,4,6,23,142,59,95,330,277)
\end{tabular}

\begin{tabular}{@{}llllll}
{\bf 38} & (63,63,95,175,87,57,24,0,0,2,0) &
(24,4) & (90,82,47,74,55,70,43,48,57,44,53,29)  &
73 & (0,1,2,3,4,11,30,53,58,6,14,36) \\
145 & (0,1,2,3,4,8,68,11,124,19,115,103)  &
217 & (0,1,2,3,4,6,79,184,98,31,48,70)  &
433 & (0,1,2,3,4,6,42,354,263,126,275,97)
\end{tabular}

\begin{tabular}{@{}llllll}
{\bf 39} & (63,63,95,175,103,27,20,8,0,0,0) &
(24,4) & (2,30,15,62,95,66,77,68,85,84,5,61)  &
73 & (0,1,2,3,4,13,36,5,15,30,20,52) \\
145 & (0,1,2,3,4,10,68,58,47,86,113,125)  &
217 & (0,1,2,3,4,6,65,92,36,180,15,144)  &
433 & (0,1,2,3,4,6,10,405,140,333,69,216)
\end{tabular}

\begin{tabular}{@{}llllll}
{\bf 40} & (63,63,95,175,103,51,16,8,4,0,0) &
(24,4) & (95,55,35,54,75,42,89,12,53,80,64,92)  &
73 & (0,1,2,3,4,15,58,46,66,52,41,71) \\
145 & (0,1,2,3,4,6,131,53,83,29,42,140)  &
217 & (0,1,2,3,4,6,21,43,150,145,99,172)  &
433 & (0,1,2,3,4,6,13,214,368,404,195,299)
\end{tabular}

\begin{tabular}{@{}llllll}
{\bf 41} & (63,63,95,175,103,57,2,12,0,0,0) &
(24,4) & (48,36,84,76,59,82,89,51,57,22,94,37)  &
73 & (0,1,2,3,4,10,35,64,71,24,17,47) \\
145 & (0,1,2,3,4,6,23,75,138,113,102,91)  &
217 & (0,1,2,3,4,6,15,98,37,197,75,49)  &
433 & (0,1,2,3,4,6,13,309,382,174,160,330)
\end{tabular}

\begin{tabular}{@{}llllll}
{\bf 42} & (63,63,95,175,103,57,8,6,0,0,0) &
(24,4) & (21,81,37,1,84,19,18,66,12,23,47,42)  &
73 & (0,1,2,3,4,18,12,35,5,45,50,55) \\
145 & (0,1,2,3,4,7,139,78,90,35,48,24)  &
217 & (0,1,2,3,4,6,24,162,139,32,208,179)  &
433 & (0,1,2,3,4,6,16,113,69,192,216,172)
\end{tabular}

\begin{tabular}{@{}llllll}
{\bf 43} & (63,63,95,175,103,57,16,2,4,0,0) &
(24,4) & (36,93,16,88,26,85,3,18,46,71,23,35)  &
73 & (0,1,2,3,4,11,25,5,13,36,31,19) \\
145 & (0,1,2,3,4,6,61,78,137,102,22,115)  &
217 & (0,1,2,3,4,6,21,69,177,147,96,141)  &
433 & (0,1,2,3,4,6,11,49,353,193,278,265)
\end{tabular}

\begin{tabular}{@{}llllll}
{\bf 44} & (63,63,95,175,103,57,16,8,0,2,0) &
(24,4) & (50,42,9,70,85,78,11,8,59,7,80,91)  &
73 & (0,1,2,3,4,11,25,5,13,19,42,47) \\
145 & (0,1,2,3,4,6,131,53,83,29,42,9)  &
217 & (0,1,2,3,4,6,95,41,17,33,154,141)  &
433 & (0,1,2,3,4,6,29,300,208,46,335,170)
\end{tabular}

\begin{tabular}{@{}llllll}
{\bf 45} & (63,63,95,175,115,51,4,8,4,0,0) &
(24,4) & (84,56,31,4,3,28,41,10,45,14,33,37)  &
73 & (0,1,2,3,4,19,49,61,55,5,41,66) \\
145 & (0,1,2,3,4,10,122,142,105,55,117,130)  &
217 & (0,1,2,3,4,6,162,186,20,108,80,136)  &
433 & (0,1,2,3,4,6,13,298,133,173,32,106)
\end{tabular}

\begin{tabular}{@{}llllll}
{\bf 46} & (63,63,95,175,115,53,2,8,4,0,0) &
(24,4) & (33,49,40,60,73,87,78,14,39,6,71,90)  &
73 & (0,1,2,3,4,6,9,26,40,13,32,58) \\
145 & (0,1,2,3,4,6,120,61,53,132,77,40)  &
217 & (0,1,2,3,4,6,33,15,178,141,67,128)  &
433 & (0,1,2,3,4,6,11,225,403,267,193,69)
\end{tabular}

\begin{tabular}{@{}llllll}
{\bf 47} & (63,63,95,175,115,53,4,8,2,0,0) &
(24,4) & (59,15,22,43,61,83,68,20,9,69,88,80)  &
73 & (0,1,2,3,4,6,59,42,63,25,47,9) \\
145 & (0,1,2,3,4,7,27,114,71,100,138,90)  &
217 & (0,1,2,3,4,6,21,122,71,63,76,49)  &
433 & (0,1,2,3,4,6,11,49,276,83,314,193)
\end{tabular}

\begin{tabular}{@{}llllll}
{\bf 48} & (63,63,95,175,115,53,8,8,0,2,0) &
(24,4) & (24,76,72,87,48,63,13,1,90,18,77,89)  &
73 & (0,1,2,3,4,10,48,70,34,12,58,53) \\
145 & (0,1,2,3,4,6,105,33,16,126,71,9)  &
217 & (0,1,2,3,4,6,15,127,66,144,42,34)  &
433 & (0,1,2,3,4,6,11,84,170,214,244,394)
\end{tabular}

\begin{tabular}{@{}llllll}
{\bf 49} & (63,63,95,175,115,53,16,0,6,0,0) &
(24,4) & (64,48,69,44,85,84,27,50,22,51,75,39)  &
73 & (0,1,2,3,4,6,31,41,14,9,58,53) \\
145 & (0,1,2,3,4,6,71,37,89,140,47,126)  &
217 & (0,1,2,3,4,6,64,37,203,108,89,141)  &
433 & (0,1,2,3,4,6,11,221,328,265,243,394)
\end{tabular}

\begin{tabular}{@{}llllll}
{\bf 50} & (63,63,95,175,115,60,1,2,4,0,0) &
(24,4) & (75,35,51,79,70,65,68,52,4,94,25,13)  &
73 & (0,1,2,3,4,6,53,41,9,31,14,58) \\
145 & (0,1,2,3,4,6,116,124,9,91,46,129)  &
217 & (0,1,2,3,4,6,33,15,41,95,141,67)  &
433 & (0,1,2,3,4,6,9,206,34,44,19,174)
\end{tabular}

\begin{tabular}{@{}llllll}
{\bf 51} & (63,63,95,175,115,60,1,4,2,0,0) &
(24,4) & (43,23,39,47,72,18,49,53,40,80,90,34)  &
73 & (0,1,2,3,4,9,64,71,17,28,53,35) \\
145 & (0,1,2,3,4,6,110,29,82,129,95,140)  &
217 & (0,1,2,3,4,6,34,21,98,48,76,71)  &
433 & (0,1,2,3,4,6,9,197,419,367,253,216)
\end{tabular}

\begin{tabular}{@{}llllll}
{\bf 52} & (63,63,95,175,115,60,4,1,2,0,0) &
(24,4) & (31,55,75,71,25,0,26,78,61,33,44,16)  &
73 & (0,1,2,3,4,10,35,70,59,12,22,30) \\
145 & (0,1,2,3,4,6,90,78,11,47,126,114)  &
217 & (0,1,2,3,4,6,84,212,164,190,40,199)  &
433 & (0,1,2,3,4,6,11,69,267,216,191,331)
\end{tabular}

\begin{tabular}{@{}llllll}
{\bf 53} & (63,63,95,175,115,60,16,0,1,2,0) &
(24,4) & (90,2,47,74,21,78,15,24,29,73,88,40)  &
73 & (0,1,2,3,4,15,22,65,35,59,28,71) \\
145 & (0,1,2,3,4,6,84,120,17,52,108,77)  &
217 & (0,1,2,3,4,6,57,35,14,179,130,153)  &
433 & (0,1,2,3,4,6,11,49,265,353,193,159)
\end{tabular}

\begin{tabular}{@{}llllll}
{\bf 55} & (63,63,95,207,87,35,28,8,0,0,0) &
(24,4) & (73,65,2,81,58,19,18,83,76,4,28,52)  &
73 & (0,1,2,3,4,6,11,35,67,47,54,61) \\
145 & (0,1,2,3,4,6,26,107,93,130,117,140)  &
217 & (0,1,2,3,4,6,23,91,206,67,41,190)  &
433 & (0,1,2,3,4,6,13,33,245,392,208,266)
\end{tabular}

\begin{tabular}{@{}llllll}
{\bf 56} & (63,63,95,207,87,41,26,4,0,0,0) &
(24,4) & (6,78,9,84,81,86,13,68,35,67,63,15)  &
73 & (0,1,2,3,4,6,13,40,65,45,19,26) \\
145 & (0,1,2,3,4,6,65,70,47,126,138,35)  &
217 & (0,1,2,3,4,6,23,65,92,180,189,16)  &
433 & (0,1,2,3,4,6,13,73,109,316,176,304)
\end{tabular}

\begin{tabular}{@{}llllll}
{\bf 57} & (63,63,95,207,87,56,25,0,0,2,0) &
(24,4) & (5,77,49,29,32,47,56,22,82,71,78,18)  &
73 & (0,1,2,3,4,6,48,41,30,12,17,54) \\
145 & (0,1,2,3,4,8,83,22,137,36,103,59)  &
217 & (0,1,2,3,4,6,41,24,108,32,162,84)  &
433 & (0,1,2,3,4,6,49,197,353,62,263,278)
\end{tabular}

\begin{tabular}{@{}llllll}
{\bf 58} & (63,63,95,207,103,19,28,8,0,0,0) &
(24,4) & (9,5,24,60,79,89,20,50,27,58,10,18)  &
73 & (0,1,2,3,4,10,50,35,71,64,32,57) \\
145 & (0,1,2,3,4,8,80,83,107,138,54,131)  &
217 & (0,1,2,3,4,6,93,23,159,122,172,184)  &
433 & (0,1,2,3,4,6,12,207,31,48,158,316)
\end{tabular}

\begin{tabular}{@{}llllll}
{\bf 60} & (63,63,95,207,103,25,26,4,0,0,0) &
(24,4) & (23,79,63,87,74,16,89,58,34,60,0,17)  &
73 & (0,1,2,3,4,6,49,39,32,16,22,68) \\
145 & (0,1,2,3,4,6,15,93,105,34,10,22)  &
217 & (0,1,2,3,4,6,10,37,98,48,199,147)  &
433 & (0,1,2,3,4,6,12,393,31,155,209,316)
\end{tabular}

\begin{tabular}{@{}llllll}
{\bf 61} & (63,63,95,207,103,35,24,12,0,0,0) &
(24,4) & (17,57,26,54,73,11,46,64,67,31,48,56)  &
73 & (0,1,2,3,4,6,38,24,8,64,21,48) \\
145 & (0,1,2,3,4,6,56,7,70,120,19,108)  &
217 & (0,1,2,3,4,6,91,200,84,190,152,164)  &
433 & (0,1,2,3,4,6,60,238,413,127,158,341)
\end{tabular}

\begin{tabular}{@{}llllll}
{\bf 62} & (63,63,95,207,103,41,18,12,0,0,0) &
(24,4) & (53,77,48,24,43,37,8,55,78,38,86,74)  &
73 & (0,1,2,3,4,6,27,7,31,58,36,14) \\
145 & (0,1,2,3,4,6,45,77,130,66,140,29)  &
217 & (0,1,2,3,4,6,42,170,24,133,13,65)  &
433 & (0,1,2,3,4,6,12,264,96,155,310,249)
\end{tabular}

\begin{tabular}{@{}llllll}
{\bf 63} & (63,63,95,207,103,41,24,6,0,0,0) &
(24,4) & (62,58,41,73,84,10,1,92,56,15,55,83)  &
73 & (0,1,2,3,4,25,38,26,65,6,45,18) \\
145 & (0,1,2,3,4,6,60,16,138,22,106,117)  &
217 & (0,1,2,3,4,6,211,103,17,130,180,155)  &
433 & (0,1,2,3,4,6,12,168,69,97,404,264)
\end{tabular}

\begin{tabular}{@{}llllll}
{\bf 64} & (63,63,95,207,103,49,10,12,0,0,0) &
(24,4) & (37,29,12,8,77,67,82,86,75,11,91,6)  &
73 & (0,1,2,3,4,6,22,26,45,36,68,62) \\
145 & (0,1,2,3,4,6,14,22,102,78,87,114)  &
217 & (0,1,2,3,4,6,23,35,94,99,151,172)  &
433 & (0,1,2,3,4,6,40,97,109,304,374,263)
\end{tabular}

\begin{tabular}{@{}llllll}
{\bf 65} & (63,63,95,207,103,49,18,8,4,0,0) &
(24,4) & (22,70,62,39,78,8,45,27,24,69,56,49)  &
73 & (0,1,2,3,4,21,61,36,5,11,47,53) \\
145 & (0,1,2,3,4,6,41,27,66,77,87,53)  &
217 & (0,1,2,3,4,6,23,30,92,207,67,40)  &
433 & (0,1,2,3,4,6,15,321,357,138,43,227)
\end{tabular}

\begin{tabular}{@{}llllll}
{\bf 66} & (63,63,95,207,103,49,24,2,4,0,0) &
(24,4) & (13,25,54,14,33,7,50,88,55,64,27,36)  &
73 & (0,1,2,3,4,10,40,24,71,13,56,48) \\
145 & (0,1,2,3,4,6,33,86,137,102,23,72)  &
217 & (0,1,2,3,4,6,17,29,39,154,176,130)  &
433 & (0,1,2,3,4,6,21,73,315,326,127,188)
\end{tabular}

\begin{tabular}{@{}llllll}
{\bf 67} & (63,63,95,207,103,49,24,8,0,2,0) &
(24,4) & (27,11,35,86,91,33,24,66,17,61,72,12)  &
73 & (0,1,2,3,4,16,40,58,65,6,47,54) \\
145 & (0,1,2,3,4,9,66,130,77,29,52,40)  &
217 & (0,1,2,3,4,6,92,195,181,48,64,13)  &
433 & (0,1,2,3,4,6,19,172,214,54,208,32)
\end{tabular}

\begin{tabular}{@{}llllll}
{\bf 68} & (63,63,95,207,103,56,3,12,0,0,0) &
(24,4) & (49,41,61,60,79,50,80,0,23,54,58,94)  &
73 & (0,1,2,3,4,16,38,70,18,50,44,66) \\
145 & (0,1,2,3,4,6,27,136,52,89,66,72)  &
217 & (0,1,2,3,4,6,45,29,90,155,183,208)  &
433 & (0,1,2,3,4,6,14,244,309,84,158,388)
\end{tabular}

\begin{tabular}{@{}llllll}
{\bf 69} & (63,63,95,207,103,56,9,6,0,0,0) &
(24,4) & (93,25,90,42,48,79,66,76,32,47,67,27)  &
73 & (0,1,2,3,4,17,53,19,65,35,70,30) \\
145 & (0,1,2,3,4,6,104,128,66,95,139,119)  &
217 & (0,1,2,3,4,6,34,74,212,98,198,48)  &
433 & (0,1,2,3,4,6,11,136,154,247,368,102)
\end{tabular}

\begin{tabular}{@{}llllll}
{\bf 70} & (63,63,95,207,103,56,17,2,4,0,0) &
(24,4) & (81,73,72,4,90,83,60,38,18,23,75,3)  &
73 & (0,1,2,3,4,11,38,5,13,32,26,18) \\
145 & (0,1,2,3,4,6,60,71,35,139,102,126)  &
217 & (0,1,2,3,4,6,14,30,164,42,94,152)  &
433 & (0,1,2,3,4,6,11,49,245,265,159,82)
\end{tabular}

\begin{tabular}{@{}llllll}
{\bf 72} & (63,63,95,207,103,56,24,1,0,2,0) &
(24,4) & (57,53,89,85,86,48,3,50,74,24,67,80)  &
73 & (0,1,2,3,4,12,49,14,23,7,34,58) \\
145 & (0,1,2,3,4,8,88,10,71,23,36,102)  &
217 & (0,1,2,3,4,6,80,208,108,38,133,157)  &
433 & (0,1,2,3,4,6,32,83,263,11,158,311)
\end{tabular}

\begin{tabular}{@{}llllll}
{\bf 73} & (63,63,95,207,115,13,22,8,0,0,0) &
(24,4) & (7,83,3,75,90,53,0,65,56,42,73,8)  &
73 & (0,1,2,3,4,6,13,23,59,68,42,49) \\
145 & (0,1,2,3,4,6,21,10,63,35,47,71)  &
217 & (0,1,2,3,4,6,103,204,128,34,154,41)  &
433 & (0,1,2,3,4,6,12,95,126,186,53,160)
\end{tabular}

\begin{tabular}{@{}llllll}
{\bf 74} & (63,63,95,207,115,19,28,0,4,0,0) &
(24,4) & (46,38,53,85,70,72,29,95,64,91,39,75)  &
73 & (0,1,2,3,4,21,35,61,10,46,6,53) \\
145 & (0,1,2,3,4,7,32,125,132,138,71,90)  &
217 & (0,1,2,3,4,6,94,212,46,122,74,199)  &
433 & (0,1,2,3,4,6,17,97,42,176,295,374)
\end{tabular}

\begin{tabular}{@{}llllll}
{\bf 75} & (63,63,95,207,115,21,14,8,0,0,0) &
(24,4) & (94,22,53,77,54,60,49,11,12,31,19,79)  &
73 & (0,1,2,3,4,6,38,51,67,33,63,21) \\
145 & (0,1,2,3,4,6,10,61,77,102,23,35)  &
217 & (0,1,2,3,4,6,13,94,196,146,69,44)  &
433 & (0,1,2,3,4,6,42,214,207,34,139,315)
\end{tabular}

\begin{tabular}{@{}llllll}
{\bf 76} & (63,63,95,207,115,21,26,0,4,0,0) &
(24,4) & (70,90,26,85,22,56,47,5,84,75,44,83)  &
73 & (0,1,2,3,4,6,52,15,68,56,42,30) \\
145 & (0,1,2,3,4,6,15,96,62,120,108,72)  &
217 & (0,1,2,3,4,6,12,98,73,24,188,48)  &
433 & (0,1,2,3,4,6,12,32,209,96,156,172)
\end{tabular}

\begin{tabular}{@{}llllll}
{\bf 77} & (63,63,95,207,115,21,28,0,2,0,0) &
(24,4) & (77,44,36,80,65,62,22,11,71,43,94,75)  &
73 & (0,1,2,3,4,10,56,24,71,40,29,64) \\
145 & (0,1,2,3,4,6,27,11,105,55,33,21)  &
217 & (0,1,2,3,4,6,13,154,176,92,190,40)  &
433 & (0,1,2,3,4,6,22,208,266,352,304,403)
\end{tabular}

\begin{tabular}{@{}llllll}
{\bf 79} & (63,63,95,207,115,28,21,0,2,0,0) &
(24,4) & (71,11,27,34,55,81,21,84,16,36,93,28)  &
73 & (0,1,2,3,4,10,38,64,51,71,31,59) \\
145 & (0,1,2,3,4,6,56,70,140,119,20,41)  &
217 & (0,1,2,3,4,6,17,174,136,188,197,163)  &
433 & (0,1,2,3,4,6,12,313,137,388,208,103)
\end{tabular}

\begin{tabular}{@{}llllll}
{\bf 80} & (63,63,95,207,115,28,28,0,0,0,1) &
(24,4) & (41,29,57,38,17,27,15,24,56,80,63,76)  &
73 & (0,1,2,3,4,34,40,28,5,53,11,59) \\
145 & (0,1,2,3,4,19,135,127,57,103,6,118)  &
217 & (0,1,2,3,4,6,31,212,21,147,49,123)  &
433 & (0,1,2,3,4,6,66,401,50,135,95,185)
\end{tabular}

\begin{tabular}{@{}llllll}
{\bf 81} & (63,63,95,207,115,35,12,12,0,0,0) &
(24,4) & (66,39,15,63,48,54,21,69,4,20,25,41)  &
73 & (0,1,2,3,4,6,38,44,16,49,56,67) \\
145 & (0,1,2,3,4,6,72,130,125,137,33,105)  &
217 & (0,1,2,3,4,6,21,62,146,185,122,172)  &
433 & (0,1,2,3,4,6,48,34,114,267,87,288)
\end{tabular}

\begin{tabular}{@{}llllll}
{\bf 82} & (63,63,95,207,115,35,20,8,4,0,0) &
(24,4) & (59,15,6,58,35,5,70,57,52,73,84,92)  &
73 & (0,1,2,3,4,10,38,18,44,71,25,64) \\
145 & (0,1,2,3,4,6,59,127,132,138,34,44)  &
217 & (0,1,2,3,4,6,91,202,123,56,72,190)  &
433 & (0,1,2,3,4,6,22,170,157,404,208,352)
\end{tabular}

\begin{tabular}{@{}llllll}
{\bf 83} & (63,63,95,207,115,37,10,12,0,0,0) &
(24,4) & (43,11,64,24,45,51,88,34,81,9,74,46)  &
73 & (0,1,2,3,4,10,38,31,5,51,59,64) \\
145 & (0,1,2,3,4,6,27,33,137,86,47,72)  &
217 & (0,1,2,3,4,6,39,126,30,47,57,70)  &
433 & (0,1,2,3,4,6,14,95,50,126,85,280)
\end{tabular}

\begin{tabular}{@{}llllll}
{\bf 84} & (63,63,95,207,115,37,12,10,0,0,0) &
(24,4) & (16,8,56,44,33,95,58,30,15,77,39,46)  &
73 & (0,1,2,3,4,6,38,56,16,26,46,68) \\
145 & (0,1,2,3,4,7,31,78,125,46,59,139)  &
217 & (0,1,2,3,4,6,12,22,49,99,45,74)  &
433 & (0,1,2,3,4,6,11,188,171,334,44,141)
\end{tabular}

\begin{tabular}{@{}llllll}
{\bf 85} & (63,63,95,207,115,37,18,8,4,0,0) &
(24,4) & (81,17,73,44,77,3,86,24,75,14,63,2)  &
73 & (0,1,2,3,4,11,38,55,71,26,44,63) \\
145 & (0,1,2,3,4,6,89,11,125,35,72,102)  &
217 & (0,1,2,3,4,6,80,212,84,108,60,34)  &
433 & (0,1,2,3,4,6,12,310,96,249,43,265)
\end{tabular}

\begin{tabular}{@{}llllll}
{\bf 86} & (63,63,95,207,115,37,20,8,2,0,0) &
(24,4) & (2,18,6,14,60,23,85,63,20,44,83,13)  &
73 & (0,1,2,3,4,10,28,64,17,71,42,54) \\
145 & (0,1,2,3,4,6,45,77,82,29,138,21)  &
217 & (0,1,2,3,4,6,23,65,180,92,31,16)  &
433 & (0,1,2,3,4,6,12,310,96,42,252,266)
\end{tabular}

\begin{tabular}{@{}llllll}
{\bf 87} & (63,63,95,207,115,37,24,2,4,0,0) &
(24,4) & (60,28,21,45,40,2,5,42,6,91,15,11)  &
73 & (0,1,2,3,4,12,38,50,16,5,56,45) \\
145 & (0,1,2,3,4,6,29,70,18,40,64,96)  &
217 & (0,1,2,3,4,6,125,48,74,135,37,23)  &
433 & (0,1,2,3,4,6,14,327,279,127,85,388)
\end{tabular}

\begin{tabular}{@{}llllll}
{\bf 88} & (63,63,95,207,115,37,24,4,2,0,0) &
(24,4) & (53,57,69,36,81,35,18,39,47,62,94,38)  &
73 & (0,1,2,3,4,6,38,44,22,61,50,11) \\
145 & (0,1,2,3,4,6,15,139,93,105,33,129)  &
217 & (0,1,2,3,4,6,14,152,30,164,94,42)  &
433 & (0,1,2,3,4,6,11,51,367,354,155,317)
\end{tabular}

\begin{tabular}{@{}llllll}
{\bf 89} & (63,63,95,207,115,37,24,8,0,2,0) &
(24,4) & (22,14,42,26,48,11,49,3,56,16,39,89)  &
73 & (0,1,2,3,4,11,38,55,71,26,32,63) \\
145 & (0,1,2,3,4,6,80,64,139,41,54,29)  &
217 & (0,1,2,3,4,6,14,150,30,125,198,48)  &
433 & (0,1,2,3,4,6,12,63,326,45,310,248)
\end{tabular}

\begin{tabular}{@{}llllll}
{\bf 90} & (63,63,95,207,115,44,3,12,0,0,0) &
(24,4) & (45,6,29,34,1,44,7,40,15,47,76,80)  &
73 & (0,1,2,3,4,17,38,70,19,50,32,65) \\
145 & (0,1,2,3,4,6,124,89,52,78,18,84)  &
217 & (0,1,2,3,4,7,172,147,24,61,35,99)  &
433 & (0,1,2,3,4,6,40,262,148,343,109,374)
\end{tabular}

\begin{tabular}{@{}llllll}
{\bf 91} & (63,63,95,207,115,44,5,10,0,0,0) &
(24,4) & (16,84,0,39,76,54,11,95,13,45,49,18)  &
73 & (0,1,2,3,4,6,17,26,21,67,36,31) \\
145 & (0,1,2,3,4,6,84,31,130,11,54,106)  &
217 & (0,1,2,3,4,6,15,91,43,122,68,31)  &
433 & (0,1,2,3,4,6,19,55,50,207,277,184)
\end{tabular}

\begin{tabular}{@{}llllll}
{\bf 92} & (63,63,95,207,115,44,12,3,0,0,0) &
(24,4) & (75,23,15,66,69,60,6,54,41,89,16,64)  &
73 & (0,1,2,3,4,6,38,68,46,58,14,26) \\
145 & (0,1,2,3,4,7,95,107,65,139,18,28)  &
217 & (0,1,2,3,4,6,80,157,108,208,38,133)  &
433 & (0,1,2,3,4,6,40,101,48,208,85,313)
\end{tabular}

\begin{tabular}{@{}llllll}
{\bf 93} & (63,63,95,207,115,44,17,2,4,0,0) &
(24,4) & (1,89,5,53,24,35,58,75,88,94,26,43)  &
73 & (0,1,2,3,4,25,38,44,5,66,58,49) \\
145 & (0,1,2,3,4,11,82,56,5,30,107,132)  &
217 & (0,1,2,3,4,6,115,129,40,208,197,157)  &
433 & (0,1,2,3,4,6,19,356,172,208,32,279)
\end{tabular}

\begin{tabular}{@{}llllll}
{\bf 94} & (63,63,95,207,115,44,17,4,2,0,0) &
(24,4) & (92,84,81,60,35,42,65,70,75,39,34,94)  &
73 & (0,1,2,3,4,6,38,44,22,61,27,66) \\
145 & (0,1,2,3,4,7,107,81,9,95,19,118)  &
217 & (0,1,2,3,4,6,11,72,122,200,31,48)  &
433 & (0,1,2,3,4,6,12,60,174,96,281,310)
\end{tabular}

\begin{tabular}{@{}llllll}
{\bf 95} & (63,63,95,207,115,44,20,1,2,0,0) &
(24,4) & (45,61,29,93,48,50,31,22,32,56,94,7)  &
73 & (0,1,2,3,4,7,49,19,9,23,36,63) \\
145 & (0,1,2,3,4,6,82,119,41,140,55,129)  &
217 & (0,1,2,3,4,6,129,24,184,15,64,172)  &
433 & (0,1,2,3,4,6,11,268,68,136,101,247)
\end{tabular}

\begin{tabular}{@{}llllll}
{\bf 96} & (63,63,95,207,115,44,20,8,0,0,1) &
(24,4) & (7,95,35,63,32,69,2,57,0,20,9,18)  &
73 & (0,1,2,3,4,17,38,45,5,11,20,51) \\
145 & (0,1,2,3,4,8,104,65,19,37,25,95)  &
217 & (0,1,2,3,4,6,57,116,130,41,13,153)  &
433 & (0,1,2,3,4,6,29,195,403,368,12,299)
\end{tabular}

\begin{tabular}{@{}llllll}
{\bf 97} & (63,63,95,207,115,49,6,8,4,0,0) &
(24,4) & (34,18,85,5,72,94,49,23,24,87,31,55)  &
73 & (0,1,2,3,4,6,38,46,16,26,55,68) \\
145 & (0,1,2,3,4,6,71,126,139,113,102,60)  &
217 & (0,1,2,3,4,6,62,147,123,200,130,175)  &
433 & (0,1,2,3,4,6,11,229,247,171,43,139)
\end{tabular}

\begin{tabular}{@{}llllll}
{\bf 98} & (63,63,95,207,115,49,12,2,4,0,0) &
(24,4) & (2,26,46,13,74,52,89,57,0,75,56,71)  &
73 & (0,1,2,3,4,24,50,41,15,31,66,6) \\
145 & (0,1,2,3,4,7,47,138,72,126,114,59)  &
217 & (0,1,2,3,4,6,34,48,179,196,151,125)  &
433 & (0,1,2,3,4,6,34,308,111,173,354,290)
\end{tabular}

\begin{tabular}{@{}llllll}
{\bf 99} & (63,63,95,207,115,49,12,4,2,0,0) &
(24,4) & (46,86,58,6,39,65,80,36,59,1,77,88)  &
73 & (0,1,2,3,4,6,14,68,26,36,46,58) \\
145 & (0,1,2,3,4,6,15,102,114,22,89,72)  &
217 & (0,1,2,3,4,6,18,103,127,67,178,59)  &
433 & (0,1,2,3,4,6,16,24,133,184,311,384)
\end{tabular}

\begin{tabular}{@{}llllll}
{\bf 100} & (63,63,95,207,115,49,20,0,6,0,0) &
(24,4) & (1,81,36,20,61,3,0,10,54,31,7,22)  &
73 & (0,1,2,3,4,16,58,41,46,11,53,71) \\
145 & (0,1,2,3,4,7,51,112,65,59,78,124)  &
217 & (0,1,2,3,4,6,15,55,93,153,122,184)  &
433 & (0,1,2,3,4,6,12,171,311,394,97,187)
\end{tabular}

\begin{tabular}{@{}llllll}
{\bf 101} & (63,63,95,207,115,52,3,8,4,0,0) &
(24,4) & (72,51,64,92,43,5,22,62,69,38,21,37)  &
73 & (0,1,2,3,4,6,20,62,24,42,30,68) \\
145 & (0,1,2,3,4,6,89,71,139,23,102,114)  &
217 & (0,1,2,3,4,6,19,136,171,162,33,197)  &
433 & (0,1,2,3,4,6,12,187,125,42,97,171)
\end{tabular}

\begin{tabular}{@{}llllll}
{\bf 102} & (63,63,95,207,115,52,5,8,2,0,0) &
(24,4) & (18,74,94,93,26,4,5,69,55,91,44,59)  &
73 & (0,1,2,3,4,6,38,68,26,46,16,55) \\
145 & (0,1,2,3,4,6,15,59,8,126,47,78)  &
217 & (0,1,2,3,4,6,23,91,206,39,31,67)  &
433 & (0,1,2,3,4,6,13,53,326,251,135,127)
\end{tabular}

\begin{tabular}{@{}llllll}
{\bf 103} & (63,63,95,207,115,52,9,2,4,0,0) &
(24,4) & (52,68,40,19,60,22,29,35,38,89,25,18)  &
73 & (0,1,2,3,4,11,38,51,13,32,6,18) \\
145 & (0,1,2,3,4,6,82,34,57,17,127,140)  &
217 & (0,1,2,3,4,6,23,91,39,190,206,65)  &
433 & (0,1,2,3,4,6,16,394,240,384,32,302)
\end{tabular}

\begin{tabular}{@{}llllll}
{\bf 104} & (63,63,95,207,115,52,9,4,2,0,0) &
(24,4) & (66,10,58,1,74,0,77,49,51,31,64,15)  &
73 & (0,1,2,3,4,6,11,41,19,15,25,31) \\
145 & (0,1,2,3,4,6,66,95,130,117,107,138)  &
217 & (0,1,2,3,4,6,19,11,197,57,136,189)  &
433 & (0,1,2,3,4,6,13,45,80,208,247,317)
\end{tabular}

\begin{tabular}{@{}llllll}
{\bf 105} & (63,63,95,207,115,52,9,8,0,2,0) &
(24,4) & (24,28,11,31,80,6,51,37,70,13,14,86)  &
73 & (0,1,2,3,4,6,22,62,8,26,36,31) \\
145 & (0,1,2,3,4,6,71,108,8,132,29,54)  &
217 & (0,1,2,3,4,6,12,98,47,147,32,22)  &
433 & (0,1,2,3,4,6,14,233,266,304,49,99)
\end{tabular}

\begin{tabular}{@{}llllll}
{\bf 106} & (63,63,95,207,115,52,12,1,2,0,0) &
(24,4) & (8,72,88,11,30,85,83,7,38,46,25,1)  &
73 & (0,1,2,3,4,6,35,47,25,68,61,56) \\
145 & (0,1,2,3,4,6,56,119,137,85,131,42)  &
217 & (0,1,2,3,4,6,18,170,74,126,44,69)  &
433 & (0,1,2,3,4,6,16,24,343,127,55,184)
\end{tabular}

\begin{tabular}{@{}llllll}
{\bf 107} & (63,63,95,207,115,52,12,8,0,0,1) &
(24,4) & (25,17,13,64,45,35,72,12,70,10,34,95)  &
73 & (0,1,2,3,4,6,38,46,55,68,16,26) \\
145 & (0,1,2,3,4,7,72,31,138,53,17,43)  &
217 & (0,1,2,3,4,6,115,23,130,64,37,86)  &
433 & (0,1,2,3,4,6,66,197,156,127,83,261)
\end{tabular}

\begin{tabular}{@{}llllll}
{\bf 108} & (63,63,95,207,115,52,17,0,6,0,0) &
(24,4) & (12,80,9,37,8,10,85,71,6,95,63,55)  &
73 & (0,1,2,3,4,6,31,51,14,35,19,67) \\
145 & (0,1,2,3,4,6,89,72,125,35,11,47)  &
217 & (0,1,2,3,4,6,12,98,35,147,173,31)  &
433 & (0,1,2,3,4,6,11,174,100,327,302,384)
\end{tabular}

\begin{tabular}{@{}llllll}
{\bf 109} & (63,63,95,207,115,52,20,0,3,0,0) &
(24,4) & (14,74,13,45,48,90,37,43,28,12,35,83)  &
73 & (0,1,2,3,4,10,38,64,32,71,52,59) \\
145 & (0,1,2,3,4,6,78,62,11,23,46,115)  &
217 & (0,1,2,3,4,6,12,91,16,66,39,190)  &
433 & (0,1,2,3,4,6,83,314,18,132,251,403)
\end{tabular}

\begin{tabular}{@{}llllll}
{\bf 110} & (63,63,95,207,115,52,24,0,1,2,0) &
(24,4) & (86,66,54,69,34,12,93,83,48,67,4,43)  &
73 & (0,1,2,3,4,6,14,68,26,58,36,46) \\
145 & (0,1,2,3,4,6,33,86,137,72,126,102)  &
217 & (0,1,2,3,4,6,68,39,131,24,15,172)  &
433 & (0,1,2,3,4,6,13,321,176,54,160,253)
\end{tabular}

\begin{tabular}{@{}llllll}
{\bf 111} & (63,63,95,207,115,52,24,0,4,0,1) &
(24,4) & (30,46,74,63,26,24,57,39,8,29,20,49)  &
73 & (0,1,2,3,4,15,40,58,5,47,10,24) \\
145 & (0,1,2,3,4,7,74,90,103,35,11,125)  &
217 & (0,1,2,3,4,6,189,154,15,130,24,41)  &
433 & (0,1,2,3,4,6,40,160,50,95,81,86)
\end{tabular}

\begin{tabular}{@{}llllll}
{\bf 112} & (63,63,95,231,103,25,10,12,0,0,0) &
(24,4) & (33,1,5,65,74,11,35,86,28,15,4,52)  &
73 & (0,1,2,3,4,6,14,21,57,9,42,50) \\
145 & (0,1,2,3,4,7,28,130,51,77,41,65)  &
217 & (0,1,2,3,4,6,31,149,73,47,22,98)  &
433 & (0,1,2,3,4,6,29,202,45,209,121,356)
\end{tabular}

\begin{tabular}{@{}llllll}
{\bf 113} & (63,63,95,231,103,25,24,2,4,0,0) &
(24,4) & (91,75,95,71,24,57,41,54,56,66,1,6)  &
73 & (0,1,2,3,4,9,23,65,71,16,41,30) \\
145 & (0,1,2,3,4,7,72,26,65,132,96,41)  &
217 & (0,1,2,3,4,6,17,29,60,66,130,181)  &
433 & (0,1,2,3,4,6,16,73,136,310,149,51)
\end{tabular}

\begin{tabular}{@{}llllll}
{\bf 114} & (63,63,95,231,103,25,24,8,0,2,0) &
(24,4) & (56,80,64,88,13,42,18,25,67,11,74,79)  &
73 & (0,1,2,3,4,38,25,12,71,9,32,58) \\
145 & (0,1,2,3,4,6,29,63,108,17,77,96)  &
217 & (0,1,2,3,4,7,46,29,158,56,72,97)  &
433 & (0,1,2,3,4,6,9,78,418,44,317,249)
\end{tabular}

\begin{tabular}{@{}llllll}
{\bf 115} & (63,63,95,231,103,41,16,10,4,0,0) &
(24,4) & (9,73,64,40,57,87,66,31,74,59,6,46)  &
73 & (0,1,2,3,4,11,21,32,51,16,71,48) \\
145 & (0,1,2,3,4,6,23,58,33,140,47,71)  &
217 & (0,1,2,3,4,6,14,183,195,180,154,130)  &
433 & (0,1,2,3,4,6,19,384,219,252,303,126)
\end{tabular}

\begin{tabular}{@{}llllll}
{\bf 116} & (63,63,95,231,103,56,1,10,4,0,0) &
(24,4) & (59,67,35,74,79,0,29,26,69,56,20,60)  &
73 & (0,1,2,3,4,6,26,13,65,41,19,31) \\
145 & (0,1,2,3,4,6,11,95,117,41,18,83)  &
217 & (0,1,2,3,4,6,23,91,190,65,206,39)  &
433 & (0,1,2,3,4,6,19,221,229,32,172,394)
\end{tabular}

\begin{tabular}{@{}llllll}
{\bf 119} & (63,63,95,231,107,13,18,12,0,0,0) &
(24,4) & (23,11,71,30,35,57,20,53,12,36,69,4)  &
73 & (0,1,2,3,4,11,58,22,49,5,14,34) \\
145 & (0,1,2,3,4,6,80,45,22,34,10,89)  &
217 & (0,1,2,3,4,6,16,146,23,48,151,190)  &
433 & (0,1,2,3,4,6,10,24,267,86,133,337)
\end{tabular}

\begin{tabular}{@{}llllll}
{\bf 120} & (63,63,95,231,107,13,20,10,0,0,0) &
(24,4) & (38,18,2,13,74,36,84,29,31,83,92,95)  &
73 & (0,1,2,3,4,6,51,61,41,30,19,68) \\
145 & (0,1,2,3,4,6,23,10,77,113,102,88)  &
217 & (0,1,2,3,4,6,13,195,40,153,92,60)  &
433 & (0,1,2,3,4,6,12,351,244,158,310,176)
\end{tabular}

\begin{tabular}{@{}llllll}
{\bf 121} & (63,63,95,231,107,19,12,12,0,0,0) &
(24,4) & (8,88,95,55,24,18,37,81,34,38,33,93)  &
73 & (0,1,2,3,4,6,31,35,67,26,22,62) \\
145 & (0,1,2,3,4,7,72,100,138,115,27,90)  &
217 & (0,1,2,3,4,6,16,23,163,153,33,41)  &
433 & (0,1,2,3,4,6,15,244,111,263,127,158)
\end{tabular}

\begin{tabular}{@{}llllll}
{\bf 122} & (63,63,95,231,107,19,20,8,4,0,0) &
(24,4) & (78,17,38,94,0,85,43,56,75,48,71,79)  &
73 & (0,1,2,3,4,6,45,39,9,56,16,52) \\
145 & (0,1,2,3,4,6,64,29,107,137,94,129)  &
217 & (0,1,2,3,4,6,13,29,130,48,63,72)  &
433 & (0,1,2,3,4,6,19,214,392,229,383,405)
\end{tabular}

\begin{tabular}{@{}llllll}
{\bf 123} & (63,63,95,231,107,21,10,12,0,0,0) &
(24,4) & (39,79,15,23,56,41,22,92,81,70,29,18)  &
73 & (0,1,2,3,4,10,18,26,64,71,54,44) \\
145 & (0,1,2,3,4,6,19,81,95,140,108,120)  &
217 & (0,1,2,3,4,6,20,122,55,94,174,149)  &
433 & (0,1,2,3,4,6,32,287,393,419,240,264)
\end{tabular}

\begin{tabular}{@{}llllll}
{\bf 124} & (63,63,95,231,107,21,12,10,0,0,0) &
(24,4) & (5,81,37,49,58,11,80,40,47,62,35,64)  &
73 & (0,1,2,3,4,10,47,64,22,71,36,60) \\
145 & (0,1,2,3,4,6,30,69,94,140,41,20)  &
217 & (0,1,2,3,4,6,49,44,123,130,71,200)  &
433 & (0,1,2,3,4,6,25,308,190,155,428,84)
\end{tabular}

\begin{tabular}{@{}llllll}
{\bf 125} & (63,63,95,231,107,21,18,8,4,0,0) &
(24,4) & (3,30,70,18,23,80,48,17,32,53,29,77)  &
73 & (0,1,2,3,4,6,53,28,63,19,68,43) \\
145 & (0,1,2,3,4,6,19,95,81,140,41,31)  &
217 & (0,1,2,3,4,6,17,29,177,39,65,92)  &
433 & (0,1,2,3,4,6,19,214,392,229,32,383)
\end{tabular}

\begin{tabular}{@{}llllll}
{\bf 126} & (63,63,95,231,107,21,20,8,2,0,0) &
(24,4) & (37,17,93,26,81,19,95,62,8,0,83,32)  &
73 & (0,1,2,3,4,6,62,41,46,19,10,26) \\
145 & (0,1,2,3,4,6,29,40,71,131,9,94)  &
217 & (0,1,2,3,4,6,10,29,40,205,155,129)  &
433 & (0,1,2,3,4,6,32,238,122,387,134,155)
\end{tabular}

\begin{tabular}{@{}llllll}
{\bf 127} & (63,63,95,231,107,21,24,2,4,0,0) &
(24,4) & (20,12,28,21,60,42,86,19,45,55,66,23)  &
73 & (0,1,2,3,4,6,35,12,9,46,60,55) \\
145 & (0,1,2,3,4,6,44,21,60,34,10,71)  &
217 & (0,1,2,3,4,6,10,138,135,108,161,189)  &
433 & (0,1,2,3,4,6,18,267,332,11,216,69)
\end{tabular}

\begin{tabular}{@{}llllll}
{\bf 128} & (63,63,95,231,107,21,24,4,2,0,0) &
(24,4) & (73,17,1,77,10,19,87,70,72,40,23,12)  &
73 & (0,1,2,3,4,7,35,19,15,9,52,47) \\
145 & (0,1,2,3,4,7,19,46,70,89,125,139)  &
217 & (0,1,2,3,4,6,13,62,173,24,131,149)  &
433 & (0,1,2,3,4,6,9,168,340,368,189,353)
\end{tabular}

\begin{tabular}{@{}llllll}
{\bf 129} & (63,63,95,231,107,21,24,8,0,2,0) &
(24,4) & (92,28,20,49,24,14,42,9,71,11,50,51)  &
73 & (0,1,2,3,4,6,17,39,27,21,45,68) \\
145 & (0,1,2,3,4,6,68,53,42,20,140,118)  &
217 & (0,1,2,3,4,6,44,116,67,14,38,130)  &
433 & (0,1,2,3,4,6,16,394,147,300,32,208)
\end{tabular}

\begin{tabular}{@{}llllll}
{\bf 131} & (63,63,95,231,107,28,5,10,0,0,0) &
(24,4) & (38,14,78,54,31,32,20,95,69,81,24,33)  &
73 & (0,1,2,3,4,6,51,62,20,30,67,36) \\
145 & (0,1,2,3,4,6,53,11,84,130,43,30)  &
217 & (0,1,2,3,4,6,11,30,98,73,48,23)  &
433 & (0,1,2,3,4,6,9,108,366,246,120,47)
\end{tabular}

\begin{tabular}{@{}llllll}
{\bf 132} & (63,63,95,231,107,28,12,3,0,0,0) &
(24,4) & (30,94,18,78,95,13,1,68,52,79,5,9)  &
73 & (0,1,2,3,4,6,46,15,67,35,20,51) \\
145 & (0,1,2,3,4,6,126,77,59,113,102,140)  &
217 & (0,1,2,3,4,6,43,24,163,15,127,153)  &
433 & (0,1,2,3,4,6,13,113,393,154,190,317)
\end{tabular}

\begin{tabular}{@{}llllll}
{\bf 133} & (63,63,95,231,107,28,17,2,4,0,0) &
(24,4) & (68,92,88,76,11,9,58,91,22,83,49,17)  &
73 & (0,1,2,3,4,6,29,67,15,52,57,34) \\
145 & (0,1,2,3,4,6,77,108,61,16,30,96)  &
217 & (0,1,2,3,4,6,13,92,158,147,38,190)  &
433 & (0,1,2,3,4,6,9,196,352,368,103,175)
\end{tabular}

\begin{tabular}{@{}llllll}
{\bf 134} & (63,63,95,231,107,28,17,4,2,0,0) &
(24,4) & (9,13,77,89,92,23,27,72,94,22,19,70)  &
73 & (0,1,2,3,4,8,21,24,62,10,36,46) \\
145 & (0,1,2,3,4,6,65,114,124,101,89,71)  &
217 & (0,1,2,3,4,6,22,15,170,89,181,65)  &
433 & (0,1,2,3,4,6,11,171,337,309,247,328)
\end{tabular}

\begin{tabular}{@{}llllll}
{\bf 135} & (63,63,95,231,107,28,20,1,2,0,0) &
(24,4) & (57,49,69,33,38,4,47,22,19,46,12,16)  &
73 & (0,1,2,3,4,6,27,45,42,22,67,62) \\
145 & (0,1,2,3,4,6,89,107,69,132,119,139)  &
217 & (0,1,2,3,4,6,24,19,175,162,136,189)  &
433 & (0,1,2,3,4,6,9,316,62,23,158,240)
\end{tabular}

\begin{tabular}{@{}llllll}
{\bf 136} & (63,63,95,231,107,28,20,8,0,0,1) &
(24,4) & (15,55,7,47,60,41,29,40,82,64,10,17)  &
73 & (0,1,2,3,4,18,67,55,47,5,25,42) \\
145 & (0,1,2,3,4,6,30,128,82,41,10,95)  &
217 & (0,1,2,3,4,7,84,46,211,31,22,58)  &
433 & (0,1,2,3,4,6,22,342,355,62,159,210)
\end{tabular}

\begin{tabular}{@{}llllll}
{\bf 137} & (63,63,95,231,107,35,16,12,4,0,0) &
(24,4) & (32,72,4,36,87,62,53,43,73,78,69,85)  &
73 & (0,1,2,3,4,10,63,51,30,5,19,41) \\
145 & (0,1,2,3,4,6,23,7,81,137,107,117)  &
217 & (0,1,2,3,4,6,15,89,183,196,49,157)  &
433 & (0,1,2,3,4,6,29,14,103,300,208,267)
\end{tabular}

\begin{tabular}{@{}llllll}
{\bf 138} & (63,63,95,231,107,37,8,14,0,0,0) &
(24,4) & (83,43,7,55,16,1,46,70,69,72,65,62)  &
73 & (0,1,2,3,4,6,17,7,65,42,30,55) \\
145 & (0,1,2,3,4,6,50,15,76,40,66,60)  &
217 & (0,1,2,3,4,6,23,92,173,181,190,207)  &
433 & (0,1,2,3,4,6,16,322,314,245,221,394)
\end{tabular}

\begin{tabular}{@{}llllll}
{\bf 139} & (63,63,95,231,107,37,16,10,4,0,0) &
(24,4) & (18,10,55,71,64,6,81,12,69,84,9,29)  &
73 & (0,1,2,3,4,6,68,61,36,46,19,53) \\
145 & (0,1,2,3,4,6,20,7,71,120,94,41)  &
217 & (0,1,2,3,4,6,14,183,107,63,132,83)  &
433 & (0,1,2,3,4,6,11,393,188,353,282,316)
\end{tabular}

\begin{tabular}{@{}llllll}
{\bf 140} & (63,63,95,231,107,37,16,12,2,0,0) &
(24,4) & (1,29,0,37,28,31,90,19,30,14,75,70)  &
73 & (0,1,2,3,4,10,31,11,40,71,59,51) \\
145 & (0,1,2,3,4,6,117,15,70,107,10,129)  &
217 & (0,1,2,3,4,6,24,81,148,98,189,175)  &
433 & (0,1,2,3,4,6,12,168,33,368,341,186)
\end{tabular}

\begin{tabular}{@{}llllll}
{\bf 141} & (63,63,95,231,107,44,1,14,0,0,0) &
(24,4) & (74,6,26,87,38,28,13,56,69,45,29,12)  &
73 & (0,1,2,3,4,6,16,29,68,43,21,51) \\
145 & (0,1,2,3,4,6,72,136,47,61,35,23)  &
217 & (0,1,2,3,4,6,43,103,200,130,67,15)  &
433 & (0,1,2,3,4,6,11,352,73,189,109,121)
\end{tabular}

\begin{tabular}{@{}llllll}
{\bf 142} & (63,63,95,231,107,44,4,11,0,0,0) &
(24,4) & (93,57,45,73,60,63,2,12,14,64,34,87)  &
73 & (0,1,2,3,4,11,34,22,49,5,58,63) \\
145 & (0,1,2,3,4,6,78,15,96,120,62,108)  &
217 & (0,1,2,3,4,6,33,56,71,82,21,96)  &
433 & (0,1,2,3,4,6,19,123,308,214,55,253)
\end{tabular}

\begin{tabular}{@{}llllll}
{\bf 143} & (63,63,95,231,107,44,16,3,4,0,0) &
(24,4) & (21,84,17,63,1,44,94,47,23,30,6,46)  &
73 & (0,1,2,3,4,6,58,24,42,62,29,67) \\
145 & (0,1,2,3,4,6,108,86,78,16,31,53)  &
217 & (0,1,2,3,4,6,96,30,23,57,150,175)  &
433 & (0,1,2,3,4,6,11,328,366,99,173,191)
\end{tabular}

\begin{tabular}{@{}llllll}
{\bf 144} & (63,63,95,231,107,44,16,5,2,0,0) &
(24,4) & (87,51,19,72,27,13,22,52,34,90,61,74)  &
73 & (0,1,2,3,4,9,18,71,25,63,47,35) \\
145 & (0,1,2,3,4,6,96,118,78,16,108,86)  &
217 & (0,1,2,3,4,6,14,36,55,179,154,126)  &
433 & (0,1,2,3,4,6,11,334,214,298,123,394)
\end{tabular}

\begin{tabular}{@{}llllll}
{\bf 145} & (63,63,95,231,107,44,16,12,0,0,1) &
(24,4) & (31,11,78,26,75,45,4,25,40,8,65,76)  &
73 & (0,1,2,3,4,6,49,17,35,53,9,63) \\
145 & (0,1,2,3,4,10,46,136,58,131,22,115)  &
217 & (0,1,2,3,4,6,16,146,23,48,31,151)  &
433 & (0,1,2,3,4,6,13,386,118,279,126,297)
\end{tabular}

\begin{tabular}{@{}llllll}
{\bf 146} & (63,63,95,231,107,49,2,12,4,0,0) &
(24,4) & (24,40,64,48,73,82,75,35,18,61,10,63)  &
73 & (0,1,2,3,4,10,63,71,34,29,25,19) \\
145 & (0,1,2,3,4,6,23,75,138,113,102,91)  &
217 & (0,1,2,3,4,6,23,190,40,155,206,129)  &
433 & (0,1,2,3,4,6,46,315,401,111,248,83)
\end{tabular}

\begin{tabular}{@{}llllll}
{\bf 147} & (63,63,95,231,107,49,4,10,4,0,0) &
(24,4) & (79,88,16,64,86,19,9,54,81,22,1,65)  &
73 & (0,1,2,3,4,6,17,39,27,32,67,55) \\
145 & (0,1,2,3,4,6,43,16,22,30,94,63)  &
217 & (0,1,2,3,4,6,13,24,147,49,62,132)  &
433 & (0,1,2,3,4,6,11,245,110,221,334,172)
\end{tabular}

\begin{tabular}{@{}llllll}
{\bf 148} & (63,63,95,231,107,49,4,12,2,0,0) &
(24,4) & (19,75,63,91,84,57,94,10,5,93,30,34)  &
73 & (0,1,2,3,4,6,31,35,19,67,14,53) \\
145 & (0,1,2,3,4,6,53,75,61,40,29,83)  &
217 & (0,1,2,3,4,6,41,73,187,171,197,15)  &
433 & (0,1,2,3,4,6,29,298,233,343,11,334)
\end{tabular}

\begin{tabular}{@{}llllll}
{\bf 149} & (63,63,95,231,107,49,16,4,6,0,0) &
(24,4) & (25,69,58,3,14,77,52,67,55,24,32,64)  &
73 & (0,1,2,3,4,6,16,44,41,27,23,68) \\
145 & (0,1,2,3,4,7,30,17,78,41,84,53)  &
217 & (0,1,2,3,4,6,21,43,145,150,172,99)  &
433 & (0,1,2,3,4,6,11,51,342,68,156,86)
\end{tabular}

\begin{tabular}{@{}llllll}
{\bf 150} & (63,63,95,231,107,52,1,10,4,0,0) &
(24,4) & (8,24,36,28,45,95,10,62,15,5,23,79)  &
73 & (0,1,2,3,4,6,64,54,18,68,30,42) \\
145 & (0,1,2,3,4,6,80,53,64,41,139,30)  &
217 & (0,1,2,3,4,6,49,23,131,62,74,183)  &
433 & (0,1,2,3,4,6,15,42,73,322,289,299)
\end{tabular}

\begin{tabular}{@{}llllll}
{\bf 151} & (63,63,95,231,107,52,1,12,2,0,0) &
(24,4) & (19,59,79,95,18,89,48,28,33,62,20,41)  &
73 & (0,1,2,3,4,11,52,58,40,31,63,71) \\
145 & (0,1,2,3,4,7,65,75,42,16,121,96)  &
217 & (0,1,2,3,4,6,21,29,65,130,180,14)  &
433 & (0,1,2,3,4,6,13,309,328,382,174,160)
\end{tabular}

\begin{tabular}{@{}llllll}
{\bf 152} & (63,63,95,231,107,52,4,9,2,0,0) &
(24,4) & (46,38,6,90,83,72,13,3,52,51,61,4)  &
73 & (0,1,2,3,4,10,48,70,24,12,43,18) \\
145 & (0,1,2,3,4,6,29,77,61,40,96,83)  &
217 & (0,1,2,3,4,6,49,196,205,64,131,39)  &
433 & (0,1,2,3,4,6,11,214,240,313,103,352)
\end{tabular}

\begin{tabular}{@{}llllll}
{\bf 153} & (63,63,95,231,107,52,8,3,4,0,0) &
(24,4) & (81,40,28,6,36,33,79,2,18,47,27,87)  &
73 & (0,1,2,3,4,6,58,42,24,62,47,67) \\
145 & (0,1,2,3,4,6,71,15,104,120,53,87)  &
217 & (0,1,2,3,4,6,13,29,173,157,130,149)  &
433 & (0,1,2,3,4,6,13,113,393,154,190,317)
\end{tabular}

\begin{tabular}{@{}llllll}
{\bf 154} & (63,63,95,231,107,52,8,5,2,0,0) &
(24,4) & (56,52,36,27,16,30,25,7,95,37,77,90)  &
73 & (0,1,2,3,4,6,25,39,57,29,15,9) \\
145 & (0,1,2,3,4,6,55,75,29,131,108,64)  &
217 & (0,1,2,3,4,6,23,30,91,40,206,65)  &
433 & (0,1,2,3,4,6,12,340,43,126,298,249)
\end{tabular}

\begin{tabular}{@{}llllll}
{\bf 155} & (63,63,95,231,107,52,8,9,0,2,0) &
(24,4) & (10,55,19,42,15,68,17,13,32,9,76,24)  &
73 & (0,1,2,3,4,6,58,62,24,43,29,67) \\
145 & (0,1,2,3,4,6,108,95,69,140,130,118)  &
217 & (0,1,2,3,4,6,92,24,206,49,40,157)  &
433 & (0,1,2,3,4,6,13,69,355,322,282,247)
\end{tabular}

\begin{tabular}{@{}llllll}
{\bf 156} & (63,63,95,231,107,52,8,12,0,0,1) &
(24,4) & (6,18,10,94,13,24,71,11,52,85,16,43)  &
73 & (0,1,2,3,4,33,16,47,40,65,6,54) \\
145 & (0,1,2,3,4,8,21,55,117,122,37,44)  &
217 & (0,1,2,3,4,6,13,149,90,24,41,156)  &
433 & (0,1,2,3,4,6,11,225,327,245,44,265)
\end{tabular}

\begin{tabular}{@{}llllll}
{\bf 157} & (63,63,95,231,107,52,16,1,6,0,0) &
(24,4) & (21,81,68,7,56,73,18,75,31,50,82,90)  &
73 & (0,1,2,3,4,6,58,41,51,31,14,68) \\
145 & (0,1,2,3,4,6,60,71,23,101,139,114)  &
217 & (0,1,2,3,4,6,17,29,177,39,65,131)  &
433 & (0,1,2,3,4,6,19,55,393,222,125,185)
\end{tabular}

\begin{tabular}{@{}llllll}
{\bf 158} & (63,63,95,231,107,52,16,4,3,0,0) &
(24,4) & (54,2,6,19,50,88,13,63,55,9,21,32)  &
73 & (0,1,2,3,4,6,42,30,27,17,22,67) \\
145 & (0,1,2,3,4,6,131,53,42,83,9,29)  &
217 & (0,1,2,3,4,6,35,49,89,159,24,98)  &
433 & (0,1,2,3,4,6,11,171,210,45,191,327)
\end{tabular}

\begin{tabular}{@{}llllll}
{\bf 159} & (63,63,95,231,107,52,16,8,1,2,0) &
(24,4) & (92,32,0,8,71,14,69,79,50,85,70,5)  &
73 & (0,1,2,3,4,11,25,5,42,13,47,19) \\
145 & (0,1,2,3,4,6,108,78,131,96,63,120)  &
217 & (0,1,2,3,4,6,13,29,130,63,149,48)  &
433 & (0,1,2,3,4,6,16,45,366,356,316,190)
\end{tabular}

\begin{tabular}{@{}llllll}
{\bf 160} & (63,63,95,231,107,52,16,8,4,0,1) &
(24,4) & (63,39,59,71,28,61,90,48,85,64,77,94)  &
73 & (0,1,2,3,4,41,51,9,6,30,16,56) \\
145 & (0,1,2,3,4,9,118,105,66,18,25,54)  &
217 & (0,1,2,3,4,6,200,108,31,150,44,138)  &
433 & (0,1,2,3,4,6,13,153,84,322,43,247)
\end{tabular}

\begin{tabular}{@{}llllll}
{\bf 161} & (63,63,95,231,121,7,10,12,0,0,0) &
(24,4) & (45,20,76,44,41,70,66,22,31,7,27,75)  &
73 & (0,1,2,3,4,6,29,13,67,59,24,43) \\
145 & (0,1,2,3,4,7,22,122,55,139,117,105)  &
217 & (0,1,2,3,4,6,39,31,123,159,23,49)  &
433 & (0,1,2,3,4,6,25,308,135,108,98,62)
\end{tabular}

\begin{tabular}{@{}llllll}
{\bf 162} & (63,63,95,231,121,7,24,2,4,0,0) &
(24,4) & (6,42,70,78,80,73,15,25,41,20,36,55)  &
73 & (0,1,2,3,4,11,31,14,20,5,40,52) \\
145 & (0,1,2,3,4,6,80,45,86,118,137,22)  &
217 & (0,1,2,3,4,6,59,94,200,152,190,85)  &
433 & (0,1,2,3,4,6,15,120,220,143,42,252)
\end{tabular}

\begin{tabular}{@{}llllll}
{\bf 164} & (63,63,95,231,121,11,6,12,0,0,0) &
(24,4) & (10,90,50,34,60,29,71,9,63,32,4,43)  &
73 & (0,1,2,3,4,16,30,67,35,58,5,53) \\
145 & (0,1,2,3,4,6,60,81,107,10,132,117)  &
217 & (0,1,2,3,4,6,17,55,120,94,151,176)  &
433 & (0,1,2,3,4,6,17,293,136,333,377,388)
\end{tabular}

\begin{tabular}{@{}llllll}
{\bf 165} & (63,63,95,231,121,11,12,6,0,0,0) &
(24,4) & (30,42,94,77,10,20,44,88,95,7,47,23)  &
73 & (0,1,2,3,4,6,32,27,22,8,36,62) \\
145 & (0,1,2,3,4,6,77,94,85,107,131,120)  &
217 & (0,1,2,3,4,6,31,16,98,74,48,23)  &
433 & (0,1,2,3,4,6,11,85,174,53,147,326)
\end{tabular}

\begin{tabular}{@{}llllll}
{\bf 166} & (63,63,95,231,121,11,20,2,4,0,0) &
(24,4) & (44,4,72,5,56,30,63,54,77,66,71,87)  &
73 & (0,1,2,3,4,6,14,19,53,30,67,42) \\
145 & (0,1,2,3,4,6,22,95,86,30,16,43)  &
217 & (0,1,2,3,4,6,22,163,179,30,153,96)  &
433 & (0,1,2,3,4,6,9,48,172,191,229,277)
\end{tabular}

\begin{tabular}{@{}llllll}
{\bf 167} & (63,63,95,231,121,11,20,4,2,0,0) &
(24,4) & (58,14,74,28,62,43,51,16,88,57,61,53)  &
73 & (0,1,2,3,4,6,59,40,31,68,63,53) \\
145 & (0,1,2,3,4,6,23,81,105,115,130,140)  &
217 & (0,1,2,3,4,6,115,81,15,157,39,89)  &
433 & (0,1,2,3,4,6,19,170,405,279,208,355)
\end{tabular}

\begin{tabular}{@{}llllll}
{\bf 168} & (63,63,95,231,121,11,20,8,0,2,0) &
(24,4) & (73,65,1,53,18,3,84,11,40,79,52,56)  &
73 & (0,1,2,3,4,6,9,30,55,59,13,42) \\
145 & (0,1,2,3,4,7,95,82,30,129,9,42)  &
217 & (0,1,2,3,4,6,60,122,49,212,24,72)  &
433 & (0,1,2,3,4,6,11,43,244,291,353,263)
\end{tabular}

\begin{tabular}{@{}llllll}
{\bf 169} & (63,63,95,231,121,14,3,12,0,0,0) &
(24,4) & (28,52,56,20,33,54,3,78,51,87,13,14)  &
73 & (0,1,2,3,4,6,10,19,41,14,26,31) \\
145 & (0,1,2,3,4,6,138,90,64,113,101,72)  &
217 & (0,1,2,3,4,6,147,200,69,82,175,96)  &
433 & (0,1,2,3,4,6,13,315,279,309,356,389)
\end{tabular}

\begin{tabular}{@{}llllll}
{\bf 170} & (63,63,95,231,121,14,6,9,0,0,0) &
(24,4) & (27,55,59,4,11,65,37,93,68,6,54,78)  &
73 & (0,1,2,3,4,11,62,56,24,5,36,46) \\
145 & (0,1,2,3,4,7,68,122,105,22,139,33)  &
217 & (0,1,2,3,4,6,17,116,172,89,181,159)  &
433 & (0,1,2,3,4,6,12,172,308,248,373,326)
\end{tabular}

\begin{tabular}{@{}llllll}
{\bf 171} & (63,63,95,231,121,14,9,6,0,0,0) &
(24,4) & (89,25,65,70,45,0,7,72,46,59,67,80)  &
73 & (0,1,2,3,4,6,22,40,17,26,67,32) \\
145 & (0,1,2,3,4,7,54,83,30,18,41,10)  &
217 & (0,1,2,3,4,7,127,99,150,20,169,46)  &
433 & (0,1,2,3,4,6,10,73,118,44,355,247)
\end{tabular}

\begin{tabular}{@{}llllll}
{\bf 172} & (63,63,95,231,121,14,10,5,0,0,0) &
(24,4) & (31,23,55,50,60,25,49,26,62,4,8,17)  &
73 & (0,1,2,3,4,6,27,36,17,31,67,22) \\
145 & (0,1,2,3,4,7,10,83,125,94,106,18)  &
217 & (0,1,2,3,4,6,31,22,74,98,48,61)  &
433 & (0,1,2,3,4,6,15,298,133,207,386,184)
\end{tabular}

\begin{tabular}{@{}llllll}
{\bf 173} & (63,63,95,231,121,14,17,2,4,0,0) &
(24,4) & (74,50,14,52,2,3,33,68,16,89,77,7)  &
73 & (0,1,2,3,4,6,13,59,49,68,23,42) \\
145 & (0,1,2,3,4,7,119,138,131,106,66,95)  &
217 & (0,1,2,3,4,6,11,32,24,98,73,48)  &
433 & (0,1,2,3,4,6,13,73,302,32,263,340)
\end{tabular}

\begin{tabular}{@{}llllll}
{\bf 174} & (63,63,95,231,121,14,17,8,0,2,0) &
(24,4) & (2,26,62,90,91,64,81,84,41,13,31,56)  &
73 & (0,1,2,3,4,6,57,52,41,30,61,10) \\
145 & (0,1,2,3,4,7,122,56,106,32,139,19)  &
217 & (0,1,2,3,4,6,13,183,107,57,87,62)  &
433 & (0,1,2,3,4,6,33,314,393,19,252,228)
\end{tabular}

\begin{tabular}{@{}llllll}
{\bf 175} & (63,63,95,231,121,14,18,1,4,0,0) &
(24,4) & (14,94,6,58,7,84,48,9,13,67,35,72)  &
73 & (0,1,2,3,4,6,19,24,41,62,48,67) \\
145 & (0,1,2,3,4,6,23,92,140,130,102,34)  &
217 & (0,1,2,3,4,6,64,93,24,206,49,181)  &
433 & (0,1,2,3,4,6,17,84,386,309,62,288)
\end{tabular}%




\begin{thebibliography}{99}

\bibitem{AbelColbournDinitz2007} R. J. R. Abel, C. J. Colbourn and J. H. Dinitz,
    Mutually Orthogonal Latin Squares (MOLS),
    \textit{Handbook of Combinatorial Designs}, second edition (C. J. Colbourn and J. H. Dinitz, eds),
    Chapman \& Hall/CRC Press, Boca Raton, 2007, 160--193.

\bibitem{AbelGeGreigLing2009} R. J. R. Abel, G. Ge, M. Greig and A. C. H. Ling,
    Further results on $(v, \{5, w^*\}, 1)$-PBDs,
    \textit{Discrete Math.} \textbf{309} (2009), 2323--2339.

\bibitem{AdamsBillingtonRodger1994} P. Adams, E. J. Billington and C. A. Rodger,
  Pasch decompositions of lambda-fold triple systems,
  \textit{J. Combin. Math. Combin. Comput.} \textbf{15} (1994), 53--63.

\bibitem{AdamsBryant1996} P. Adams and D. E. Bryant,
  Decomposing the complete graph into Platonic graphs,
  \textit{Bull. Inst. Combin. Appl.} \textbf{17} (1996), 19--26.

\bibitem{AdamsBryantBuchanan2008} P. Adams, D. E. Bryant and M. Buchanan,
  A survey on the existence of \textit{G}-designs,
  \textit{J. Combin. Des.} \textbf{16} (2008), 373--410.

\bibitem{AdamsBryantForbesGriggs2012} P. Adams, D. E. Bryant, A. D. Forbes and T. S. Griggs,
  Decomposing the complete graph into dodecahedra,
  \textit{J. Statist. Plann. Inference} \textbf{142} (2012), 1040--1046.

\bibitem{BrouwerSchrijverHanani1977} A. E. Brouwer, A. Schrijver and H. Hanani,
    Group divisible designs with block size four,
    \textit{Discrete Math.} \textbf{20} (1977), 1--10.

\bibitem{BryantZanatiGardner1994} D. E. Bryant, S. El-Zanati and R. Gardner,
  Decompositions of $K_{m,n}$ and $K_n$ into cubes,
  \textit{Australas. J. Combin.} \textbf{9} (1994), 285--290.

\bibitem{ForbesGriggs2012} A. D. Forbes and T. S. Griggs,
  Icosahedron designs,
  \textit{Australas. J. Combin.} \textbf{52} (2012), 215--228.

\bibitem{ForbesGriggs2013} A. D. Forbes and T. S. Griggs,
  Archimedean graph designs,
  \textit{Discrete Math.} \textbf{313} (2013) 1138--1149.

\bibitem{ForbesForbesGriggs2017} A. D. Forbes, T. J. Forbes and T. S. Griggs,
  Archimedean graph designs - II,
  \textit{Discrete Math.} \textbf{340} (2017) 1598--1611.

\bibitem{Ge2007} G. Ge,
    Group Divisible Designs,
    \textit{Handbook of Combinatorial Designs}, second edition (C. J. Colbourn and J. H. Dinitz, eds),
    Chapman \& Hall/CRC Press, Boca Raton, 2007, 255--260.

\bibitem{GeLing2004c} G. Ge and A. C. H. Ling,
    Some More 5-GDDs and Optimal $(v,5,1)$-Packings,
    \textit{J. Combin. Des.} \textbf{12} (2004), 132--141).

\bibitem{GeLing2005} G. Ge and A. C. H. Ling,
    Asymptotic results on the existence of 4-RGDDs and uniform 5-GDDs,
    \textit{J. Combin. Des.} \textbf{13} (2005), 222--237.

\bibitem{GeMiao2007} G. Ge and Y. Miao,
    PBDs, Frames and Resolvability,
    \textit{Handbook of Combinatorial Designs}, second edition (C. J. Colbourn and J. H. Dinitz, eds),
    Chapman \& Hall/CRC Press, Boca Raton, 2007, 261--265.

\bibitem{GriggsResminiRosa1992} T. S. Griggs, M. J. deResmini and A. Rosa,
  Decomposing Steiner triple systems into four-line configurations,
  \textit{Ann. Discrete Math.} \textbf{52} (1992), 215--226.

\bibitem{Hanani1961} H. Hanani,
  The existence and contruction of balanced incomplete block designs,
  \textit{Ann. Math. Statist.} \textbf{32} (1961), 361--386.

\bibitem{HananiRayChaudhuriWilson1972} H. Hanani, D. K. Ray-Chaudhuri and R. M. Wilson,
  On resolvable designs,
  \textit{Discrete Math.} \textbf{3} (1972), 343--357.

\bibitem{Kirkman1847} T. P. Kirkman,
    On a Problem in Combinatorics,
    \textit{Cambridge Dublin Math. J.} \textbf{2} (1847), 191--204.

\bibitem{Kotzig1981} A. Kotzig,
  Decompositions of complete graphs into isomorphic cubes,
  \textit{J. Combin. Theory B} \textbf{31} (1981), 292--296.

\bibitem{Maheo1980} M. Maheo,
  Strongly graceful graphs,
  \textit{Discrete Math.} \textbf{29} (1980), 39--46.

\bibitem{McKay1979} B. D. McKay,
    Transitive Graphs With Fewer Than Twenty Vertices,
    \textit{Math. Comp.} \textbf{33} (1979), 1101--1121.

\bibitem{Meringer1999} M. Meringer,
    Fast Generation of Regular Graphs and Construction of Cages,
    \textit{J. Graph Theory} \textbf{30} (1999), 137--146.

\bibitem{MeringerReggraphs1999} M. Meringer,
   Regular Graphs, \\
   \url{https://www.mathe2.uni-bayreuth.de/markus/reggraphs.html#CRG}.

\bibitem{Rosa1967} A. Rosa,
   On certain valuations of the vertices of a graph,
   \textit{Theory of Graphs (Internat. Sympos., Rome, 1966)} (P. Rosenstiehl, editor),
   Gordon and Breach, New York, 1967, 349-–355.

\bibitem{WeiGe2014} H. Wei and G. Ge,
    Some more 5-GDDs, 4-frames and 4-RGDDs,
    \textit{Discrete Mathematics} \textbf{336} (2014), 7--21.

\bibitem{Wilson1972} R. M. Wilson,
    An existence theory for pairwise balanced designs I. Composition theorems and morphisms,
    \textit{J. Combin. Theory A} \textbf{13} (1972), 20--236.

\bibitem{Wilson1976} R. M. Wilson,
    Decompositions of complete graphs into subgraphs isomorphic to a given graph,
    \textit{Congr. Numer.} \textbf{15} (1976), 647--659.

\end{thebibliography}
\end{document}